**DAHLGREN DIVISION**
**NAVAL SURFACE WARFARE CENTER**
Dahlgren, Virginia 22448-5100

# STATE SPACE CONSISTENCY AND DIFFERENTIABILITY CONDITIONS FOR A CLASS OF CAUSAL DYNAMICAL INPUT-OUTPUT SYSTEMS

**BY DEMETRIOS SERAKOS**

**WARFARE SYSTEMS DEPARTMENT**

| REPORT DOCUMENTATION PAGE | | | Form Approved OMB No. 0704-0188 |
|---|---|---|---|

| 1. AGENCY USE ONLY (Leave blank) | 2. REPORT DATE | 3. REPORT TYPE AND DATES COVERED |
|---|---|---|
| | January 2008 | Final |

| 4. TITLE AND SUBTITLE | 5. FUNDING NUMBERS |
|---|---|
| STATE SPACE CONSISTENCY AND DIFFERENTIABILITY CONDITIONS FOR A CLASS OF CAUSAL DYNAMICAL INPUT-OUTPUT SYSTEMS | |
| 6. AUTHOR(s)<br>Demetrios Serakos | |
| 7. PERFORMING ORGANIZATION NAME(S) AND ADDRESS(ES)<br>NAVAL SURFACE WARFARE CENTER DAHLGREN<br>BALLISTIC MISSILE DEFENSE CODE W31<br>19008 WAYSIDE DRIVE SUITE 337<br>DAHLGREN VA 22448-5162 | 8. PERFORMING ORGANIZATION REPORT NUMBER<br>NSWCDD/TR-08/3 |
| 9. SPONSORING/MONITORING AGENCY NAME(S) AND ADDRESS(ES) | 10. SPONSORING/MONITORING AGENCY REPORT NUMBER |

11. SUPPLEMENTARY NOTES

| 12a. DISTRIBUTION/AVAILABILITY STATEMENT | 12b. DISTRIBUTION CODE |
|---|---|
| Approved for public release; distribution is unlimited. | |

13. ABSTRACT (Maximum 200 words)

A causal input-output system may be described by a function space for inputs, a function space for outputs, and a causal operator mapping the input space into the output space. A particular representation of the state of such a system at any instant has been defined as an operator from the space of possible future inputs to that of future outputs. This representation is called the natural state. The purpose of this report is to investigate additional properties of the natural state in two areas. The first area has to do with the possibility of determining the input-output system from its natural state set. A counterexample where this is not possible is given. Sufficient conditions for identifying the system from its natural state set are given. The results in this area are mostly for time-invariant systems. There are also some preliminary observations on reachability. The second area deals with differentiability properties involving the natural state inherited from the input-output system, including differentiability of the natural state and natural state trajectories. A differential equation representation is given. The results presented in this report may be considered as aids in modeling physical systems because system identification from state set holds in many models and is only tacitly assumed; also, differentiability is a useful property for many systems.

| 14. SUBJECT TERMS | 15. NUMBER OF PAGES |
|---|---|
| input-output systems, state, differential equations | 80 |
| | 16. PRICE CODE |

| 17. SECURITY CLASSIFICATION OF REPORTS | 18. SECURITY CLASSIFICATION OF THIS PAGE | 19. SECURITY CLASSIFICATION OF ABSTRACT | 20. LIMITATION OF ABSTRACT |
|---|---|---|---|
| UNCLASSIFIED | UNCLASSIFIED | UNCLASSIFIED | UU |

NSN 7540-01-280-5500

Standard Form 298 (Rev 2-89)
Prescribed by ANSI std. Z39-18
298-102

# FOREWORD

This report presents an analysis of input-output systems with regard to the mathematical concept of "state." The state of a system condenses, in a usable way, the effects of past inputs to the system. This report is the result of a theoretical examination of the state concept and the conditions under which state assumptions are valid. The results developed in this report are useful aids for understanding proper modeling of physical systems including the development of guidance and control algorithms, and in developing simulations.

The author is indebted to Professor William L. Root (1919-2007) for help with the text and with Proposition 20, and for providing Proposition 15. The author wishes to thank anonymous reviewers for various comments and for simplifying Example 8. Mr. John E. Gray of NSWCDD's Electromagnetic and Sensor Systems Department provided the first three paragraphs of the introduction.

The author wishes to thank Mr. Jaan A. Laanisto and Mr. Roger Carr for their support in getting this report published, which was a long-time goal of the author.

This report has been reviewed by John T. Beck, Head, Ballistic Missile Defense Systems Engineering Branch, and Emmanuel Skamangas, Head, Warfare Systems Development Division.

Approved by:

DONALD L. BURNETT, Head
Warfare Systems Department

# CONTENTS

## NOMENCLATURE

| | |
|---|---|
| $u, v$ | Input, Page 2-1 |
| $u_t, v_t, u^t, v^t$ | Truncated (past) input, Page 2-2 |
| $u_{t,\infty}, v_{t,\infty}$ | Future input, Page 2-3 |
| $U$ | Input space, Page 2-3 |
| $U_t$ | Truncated (past) input space, Page 2-3 |
| $U_{t,\infty}$ | Future input space, Page 2-3 |
| $Y, Y_t$ | Output space, truncated (past) output space, Page 2-4 |
| $F$ | System operator, Page 2-4 |
| $\widetilde{F}_t$ | truncated system operator, Page 2-4 |
| $F_t$ | centered truncated system operator, Page 2-4 |
| $(Y, F, U)$ | Input-output system, Definition 4 |
| $\xi, \eta$ | Natural state, Definition 7 |
| $\widetilde{\xi}_t, \widetilde{\eta}_t$ | truncated natural state, Page 2-5 |
| $\xi_t, \eta_t$ | centered truncated natural state, Page 2-5 |
| $\Xi, H$ | Natural state set, Page 2-6 |
| $\mathrm{L}_t$, $\mathrm{R}_t$ | Left-, right-shift by $t$, Definition 1 |
| $\alpha_U, K_U$ | $\alpha$, $K$ from (5) of Definition 1, associated with $U$ |
| $\mathcal{F}_N(U, Y)$ | Bounded operators $U$ to $Y$, Page 2-6 |
| $\mathcal{C}_N(U, Y)$ | Bounded, continuous operators $U$ to $Y$, Page 2-6 |
| $\lVert\cdot\rVert$ | Vector norm |
| $\lVert\cdot\rVert_{s,t}$ | Fitted Family seminorm over $(s, t)$, Definition 1 |
| $[\![ \cdot ]\!]$ | Operator norm, Equation (2-5) |
| $\leftarrow\!\!\!\!+\!\!\!\rightarrow$ | Splice, Definition 3 |
| $\square$ | Indicates the end of a proof |

# 1 INTRODUCTION

*The concept of a state in a mathematical formulation of a physical law dates back to Newton's Principia — where the first comprehensive formulation of some aspects of natural law was cast in the form of differential equations. The Hamiltonian formulation of mechanics in the early 19th century was the first explicit formulation of mechanics that could be regarded as a state space formulation with the position and momentum variables forming the two state variables. State space could not be said to be fully born until the Cartesian formulation of mechanics in terms of symplectic variables in the 30's. With the advent of the mathematical formulation of quantum mechanics and quantum field theory during the last sixty years, the concept of state has started to break down as the concept of a distinctive particle has begun to disappear to the point where many physicists are inclined to think the particle concept has outlived its usefulness. Without the concept of identity, there is a breakdown of distinctiveness that is needed to achieve a state-based formulation of dynamics.*

*The picture of what a state is in the area of mathematical engineering is considerably better. Post WWII, mathematics was used to profit in the fifties and sixties to formulate a modern mathematical theory of systems theory (culminating in Zadeh's book on systems theory [19]), the Kalman solution of the estimation problem (Kalman filter), the formulation of optimization theory, as well as Bellman's development of dynamic programming. The deeper and rigorous understanding of the state concept has made all of this progress possible.*

*As engineering extends the state concept into new areas where the concept of state has not previously been applied — it remains the task of mathematical engineering to help "clean up" the formulation of what a state is so that it can be applied to wider and wider domains. In this document, we take some concepts that were developed jointly by the author and W.L. Root to extend state concept so it may be applied to new problem domains where new state concepts are needed. At the heart of the definition of what a state is, there is an ordinary or partial differential equation that is used to formulate an explanation of the behavior of a physical system, and the behavior of certain physical systems suggests new approaches to the state concept.*[1]

State is of fundamental importance in systems theory. The exact description of state, as can be seen in the literature (some mentioned below), depends on the system. However, state condenses the past input and has a state transition property. A generic and general definition of state is given by Kalman, Falb, and Arbib, [6]. Natural state is defined for input-output systems. Natural state[2] was defined by Root, [10], and it was shown that it inherits certain continuity properties from the system. Natural state conforms to the Kalman-Falb-Arbib state.

We discussed the natural state of an input-output system in a previous paper, Root and Serakos, [13], and examples were presented. The natural state is a mapping from future inputs to future outputs. The natural state is specified by the of past state and past input. The relationship of the natural state to the Nerode equivalence classes, [7], and the more inclusive behavioral approach of Willems, [18], is discussed in [13]. The natural state inherits many properties from the input-output system. Some of these properties are discussed in [10] and [13]. For instance, the natural state is continuous (respectively bounded) if the input-output system is continuous (respectively bounded.) The mapping from past input to the natural state space is continuous under specified conditions. Conditions for the continuity of the natural state trajectory and future state, based on past state, are presented in [13]. Also in [13], the natural state of a nonlinear system, given by an integral operator, is computed, and it is shown that for linear time-invariant systems, the natural state is equivalent to the traditional state.

[1] John E. Gray, Electromagnetic and Sensor Systems Department, NSWCDD.
[2] The definition of natural state is provided in this report, Definition 7.

In this report, we discuss further properties of the natural state. Hence, in a strong sense this report is a continuation of [13]. There are two areas of emphasis: First, we define the natural state space as the set of all natural states that can be attained by an input-output system. Considering a family of input-output systems with input and output spaces in common, a map may be defined from this family to their natural state spaces. Of natural mathematical interest is the invertibility of this natural state space map. For typical systems, the natural state space map is invertible; hence, it may be a desirable property to extend to more general systems. In this light, the invertibility of this map may be important in ruling out undesirable or unintended properties in models of systems. An example is given that shows the natural state space map need not be invertible, in particular when the system has "too much memory." Essentially the same example shows that reachability of one natural state from another can be severely curtailed when there is too much memory. Finite memory or fading memory has been mentioned in connection with a number of applications, such as Kalman filtering, Anderson and Moore, [1], systems theory, Root, [9], and stability, Serakos, [15]. A proposition giving sufficient conditions for the invertibility of the natural state space map, which involves memory length, is presented. A general class of input-output systems represented by integral operators are considered to see when they are determined from their natural state space. A proposition, also involving memory, is presented that gives sufficient conditions for reachability.

The second area has to do with natural state differentiability properties. The properties are patterned after the continuity properties in [13]. The results presented are differentiability of the natural state, differentiability of the mapping from the past input to the natural state, and differentiability of the natural state trajectories. These results are inherited from appropriate differentiability properties of the input-output system. These differentiability properties may be contrasted with the work of Jakubczyk, [4], [5], and several of his references. In these papers a full state space differential equation is developed. For example, this is in contrast to the differentiability of individual natural state trajectories developed here.

The results presented in this report are original and the author's, except as stated in the foreword.

# 2 PRELIMINARIES

The purpose of this chapter is to present the reader with a workable definition of an input-output system and natural state. An input-output system is herein denoted $(Y, F, U)$ where $F$ is a mapping from an input space $U$ to an output space $Y$, and where $U$ and $Y$ are translation-invariant spaces of vector-valued time functions. The vector values of these time functions typically are $\Re^N, N = 1, 2, \cdots$, but only need to be Banach spaces. Other spaces and mappings related to $(Y, F, U)$ will be introduced, but $U$ and $Y$ with whatever affixes they carry always refer to input and output spaces. Mappings from various input spaces to output spaces are denoted either $F$ or $\xi$, again with qualifying affixes; $\xi$ (or lower case Greek) is reserved for state operators (i.e., for natural states).

The input and output space metrics are set up by seminorms referred to as fitted families (FFs) of seminorms. Roughly speaking, fitted families work like $L_p$ norms on time functions with the additional feature that a time weighting can be incorporated so that the distant past of an input or output time function may be de-emphasized. The notation $\|u_{s,t}\|_{s,t}$ indicates the norm (e.g., weighted $L_p$ norm) of the input $u$ over the interval of time $(s, t]$. $U_{s,t}$ is the space of inputs over the same interval. FFs were initially described by Root, [9].

**Definition 1 ([9])** *Let $\mathcal{L} = \mathcal{L}(\Re, E)$ be a linear space of time functions from $\Re$ into a Banach space $E$ such that any translate of a function in $\mathcal{L}$ is also a function in $\mathcal{L}$. Let $\mathcal{N} = \{\|\cdot\|_{s,t}, -\infty < s < t < \infty\}$ be a family of seminorms on $\mathcal{L}$ satisfying the following conditions:*

*(1) For $f_1, f_2 \in \mathcal{L}$, if $f_1(\tau) = f_2(\tau)$ for $s < \tau \leq t$ then $\|f_1 - f_2\|_{s,t} = 0$.*

*(2) Let $\mathrm{L}_\tau$ denote shift to the left by $\tau$. For all $f \in \mathcal{L}$, $\|\mathrm{L}_\tau f\|_{s-\tau,t-\tau} = \|f\|_{s,t}$.*

*(3) Let $r < s < t$. Then for all $f \in \mathcal{L}$, $\|f\|_{s,t} \leq \|f\|_{r,t}$.*

*(4) Let $r < s < t$. Then for all $f \in \mathcal{L}$, $\|f\|_{r,t} \leq \|f\|_{r,s} + \|f\|_{s,t}$.*

*(5) There exists $0 < \alpha \leq \infty$ and $K \geq 1$ such that if $0 < t - r \leq \alpha$ and $r < s < t$, then for all $f \in \mathcal{L}$, $\|f\|_{r,s} \leq K \|f\|_{r,t}$.*

*The pair $(\mathcal{L}, \mathcal{N})$ is called a FF of seminorms on $\mathcal{L}$. The normed linear space formed from equivalence classes of functions in $\mathcal{L}$ with norm $\|\cdot\|_{s,t}$ is denoted $A_{s,t}$. The elements of $A_{s,t}$ are the equivalence classes determined by: $f \sim g$, $f, g \in \mathcal{L}$ if and only if $\|f - g\|_{s,t} = 0$. They are denoted $u_{s,t}$, $y_{s,t}$, etc. The set $\{A_{s,t}\}$, $-\infty < s < t < \infty$, is the FF of normed linear spaces given by $(\mathcal{L}, \mathcal{N})$.*

A fairly wide class of examples of FFs are given by weighted $L_p$ spaces. For $1 \leq p < \infty$, let $w$ be a fixed nonnegative, nonincreasing Lebesgue measurable real-valued function and let $\mathcal{L} = \mathcal{L}(\Re, \Re^N)$ be the set of $N$ vector-valued functions on $\Re$ that are $p$-integrable Lebesgue on finite intervals. Then, for $f \in \mathcal{L}$ the seminorms

$$\|f\|_{s,t} = (\int_s^t \|f(\tau)\|^p \, w(t-\tau) d\tau)^{1/p} \tag{2-1}$$

satisfy Conditions (1),$\cdots$,(5) of Definition 1. For $p = \infty$, let $\mathcal{L}$ be the set of essentially bounded functions. When $\alpha = +\infty$ and $K = 1$, a uniform time weighting is given and $(\mathcal{L}, \mathcal{N})$ is referred to as a *standard* FF of seminorms. As time weighting was not essential to its purpose, standard FFs were used in [13].

We define another FF of seminorms. For $f \in \mathcal{L}$, put $\|f\|^{s,t} \triangleq \sup_{s<\tau\leq t} \|f\|_{s,\tau}$. With $\mathcal{M}$ indicating this new set of seminorms, $(\mathcal{L}, \mathcal{M})$ is indeed a FF of seminorms, [9].

A FF $(\mathcal{L}, \mathcal{N})$ and $\{A_{s,t}\}$, $-\infty < s < t < \infty$, can be augmented to include $\|\cdot\|_{-\infty,t}$ by taking the limit $s \to -\infty$, since by (3) of Definition 1 $\|f\|_{s,t}$ is monotone nondecreasing as $s \to -\infty$ with $t$ fixed. Let $\mathcal{L}_0 = \{f \in \mathcal{L} | \lim_{s\to-\infty} \|f\|_{s,t} < \infty, \; t \in \Re\}$. For $f \in \mathcal{L}_0$, define

$$\|f\|_t \triangleq \lim_{s\to-\infty} \|f\|_{s,t} \; = \|f\|_{-\infty,t} \, . \tag{2-2}$$

With the meaning of $(\mathcal{L}, \mathcal{N})$ thus extended, $\|\cdot\|_{s,t}$ is defined for $-\infty \leq s < t < \infty$. The *left-expanded FF of seminorms* is thereby defined and is denoted $(\mathcal{L}_0, \mathcal{N})$. It still satisfies all the Conditions $(1), \cdots, (5)$.

To discuss natural state we also need to define $\|\cdot\|_{s,\infty}$ and $A_{s,\infty}$ in a meaningful way. For a FF, this is done by taking the supremum. Let $\mathcal{L}_{00} = \{f \in \mathcal{L}_0 | \sup_t \|f\|_t < \infty\}$. For $f \in \mathcal{L}_{00}$ define

$$\|f\|_{s,\infty} \triangleq \sup_{t>s} \|f\|_{s,t} \;\; ; \; -\infty \leq s \, . \tag{2-3}$$

It may be readily verified that if $(\mathcal{L}, \mathcal{N})$ is a FF for indices satisfying $-\infty < s < t < \infty$ then, with definitions given by (2-2) and (2-3), $(\mathcal{L}_{00}, \mathcal{N})$ is a FF for indices satisfying $-\infty \leq s < t < \infty$ and satisfies Conditions 1, 2, 3, and 5 of Definition 1 for indices $-\infty \leq s < t \leq \infty$. (For standard FFs Condition 4 holds for both cases.) $\{(\mathcal{L}_{00}, \mathcal{N}), \; \|\cdot\|_{s,t}\,, \; -\infty \leq s < t \leq \infty\}$ is called the *expanded family of seminorms* determined by $(\mathcal{L}, \mathcal{N})$. Note that $(\mathcal{L}_{00}, \mathcal{M})$ similarly defined is a FF for $-\infty \leq s < t \leq \infty$.

For $f \in \mathcal{L}_{00}$ we put

$$\|f\| \triangleq \sup_{t\in\Re} \|f\|_t \; = \|f\|_{-\infty,\infty} \, . \tag{2-4}$$

The normed linear space consisting of equivalence classes of functions in $\mathcal{L}_{00}$ with the norm (2-4) is called the *bounding space* $A$ for the family $\{A_{s,t}\}$.

The *extended space* $A^e$ for the family $\{A_{s,t}\}$ is the set of all equivalence classes of functions $f$ in $\mathcal{L}_0$ ($f \sim g$ iff $\|f - g\| = 0$) for which $\|f\|_s < \infty$ for all $s$. It does not have a norm and, indeed, is given no topology. This definition agrees with the notion of extended space commonly used in the control literature.

It is possible that a FF $(\mathcal{L}, \mathcal{N})$ has a vacuous expansion in the sense that $\mathcal{L}_{00}$ is the empty set. An obvious example of this is given when $\mathcal{L}$ is the set of all constant real-valued functions on $\Re$ and $\mathcal{N}$ is the set of $L_1$-norms on finite intervals. To prevent this from happening and further to prevent the bounding space $A$ from being too small (in a sense to be made explicit below) we can require that a FF be "full," as in the following definition.

**Definition 2** *The FF $(\mathcal{L}, \mathcal{N})$ is full if each equivalence class $u_{s,t} \in A_{s,t}$, $-\infty < s < t < \infty$, has a representing function belonging to $\mathcal{L}_{00}$.*

When this definition is satisfied, then for all pairs $(s, t)$, $-\infty < s < t < \infty$, there is a 1:1 correspondence between the normed linear space $A_{s,t}$ determined by $(\mathcal{L}, \mathcal{N})$ and the normed linear space $A'_{s,t}$ determined by $(\mathcal{L}_{00}, \mathcal{N})$, which preserves the normed linear space structure. The correspondence is given by $u_{s,t} \leftrightarrow u'_{s,t}$, $u_{s,t} \in A_{s,t}$, $u'_{s,t} \in A'_{s,t}$ if and only if $u_{s,t}$ and $u'_{s,t}$ have a common representing function $f \in \mathcal{L}_{00}$. Thus, if $(\mathcal{L}, \mathcal{N})$ is full, we need not distinguish between $A_{s,t}$ and $A'_{s,t}$. Henceforth, every FF mentioned is assumed to be full. The FFs formed with $L_p$ spaces as described above are full. To emphasize the relations among the equivalence classes, suppose the function $f \in \mathcal{L}_{00}$ determines the equivalence classes $u \in A$, $u^t \in A^t$, $u_t \in A_t$ and $u_{s,t} \in A_{s,t}$. Since

$$\|f - g\|_{s,t} \leq \|f - g\|_t \leq \|f - g\|^t \leq \|f - g\|$$

the equivalence class $u$ considered as a set of functions is entirely contained in the equivalence class $u^t$ considered as a set of functions, and similarly $u^t \subset u_t$ and $u_t \subset u_{s,t}$; also $u \subset u^{t,\infty}$. Thus, e.g., given $t$, $u$ determines $u_t$ and $u^{t,\infty}$. Therefore, if $f$ determines $u$ and $-\infty \leq s < t \leq \infty$, it is meaningful for example to write $\|f\|_{s,t}$, $\|u\|_{s,t}$, $\|u^t\|_{s,t}$, $\|u^{s,\infty}\|_{s,t}$, $\|u_t\|_{s,t}$, $\|u_{s,t}\|_{s,t}$, and they are all equal. Let $-\infty \leq r < s < t \leq \infty$. Then, since $\|f\|_{s,t} \leq \|f\|_{r,t}$ for $f \in \mathcal{L}_{00}$, the partitioning of $\mathcal{L}_{00}$ into equivalence classes by $\|\cdot\|_{r,t}$ results in a finer partition than that given by $\|\cdot\|_{s,t}$. That is, letting $f$ determine $u \in U$, we have $u \subset u_{r,t} \subset u_{s,t}$.

In order to define the natural state, it is necessary to consider an arbitrary past input concatenated with an arbitrary future input; that is, to "splice" two inputs.

**Definition 3** *For $-\infty \leq r < s < t < \infty$, and $h$, $g \in \mathcal{L}$, the splice of $h$ and $g$ over $(r, t]$ at $s$ is defined and equals $f$ if*

$$f(\tau) = \begin{cases} h(\tau), & r < \tau \leq s \\ g(\tau), & s < \tau \leq t \end{cases}$$

*belongs to $\mathcal{L}$. It is denoted $f_{r,t} = h_{r,s} \mapsto g_{s,t}$. For $t = \infty$, the splice of $h$ and $g$ equals $f$ if*

$$f(\tau) = \begin{cases} h(\tau), & r < \tau \leq s \\ g(\tau), & s < \tau \end{cases}$$

*belongs to $\mathcal{L}$. It is denoted $f_{r,\infty} = h^{r,s} \mapsto g_{s,\infty}$.*

If $u_{r,s} \in A_{r,s}$, $v_{s,t} \in A_{s,t}$ are determined by functions $h$ and $g$ respectively, and $h_{r,s} \mapsto g_{s,t}$ exists, then for a FF, the splice of $u$ and $v$ (or $u_{r,s}$ and $v_{s,t}$) over $(r, t]$ at $s$ is defined to be the element $w_{r,t} \in A_{r,t}$ determined by $h_{r,s} \mapsto g_{s,t}$; we write $w_{r,t} = u_{r,s} \mapsto v_{s,t}$. For $t = \infty$ we write $w_{r,\infty} = u^{r,s} \mapsto v_{s,\infty}$. These are not meaningful until it is proved that the splice is independent of the particular functions $h$ and $g$ representing the equivalence classes $u_{r,s}$ and $v_{s,t}$. However, this proof follows easily from Definition 1.

Any input space $U$ is herein taken to be either the bounding space $A$ of a FF $\{A_{s,t}\}$ that permits splicing or a translation-invariant subset of $A$. The extended space $A^e$ can appear in an auxiliary role. Whether the normed linear space $A$ is complete or not is irrelevant for the purposes of this report. We write $U_{s,t}$, $-\infty \leq s < t \leq \infty$, to denote the set ("space") of equivalence classes of functions belonging to $U$ as determined by $\|\cdot\|_{s,t}$. If $U = A$, then $U_{s,t}$ is the normed linear space $A_{s,t}$; if $U$ is a subset of $A$, $U_{s,t}$ is a space only in the sense that it is a subset of $A_{s,t}$. We call any $U_{s,t}$ a *truncated input space* and write $U_t$ for $U_{-\infty,t}$.

The requirement that $A$ permits splicing means that, if $U = A$, future inputs at any $t$ can be arbitrary, independent of the past. Unfortunately, spaces of functions everywhere continuous on $\Re$ do not qualify, but this appears to be a minor drawback. It is sometimes desired that $U$ be a translation-invariant bounded (or even totally bounded) subset of $A$; we always assume $U$ contains the zero function. If $U$ is a proper subset of $A$, a splice of two elements in $U$ does not necessarily belong to $U$, of course. However, we require that for all $u \in U$, both $u^t \!\!\leftrightarrow 0_{t,\infty}$ and $0^t \!\!\leftrightarrow u_{t,\infty}$ belong to $U$.

The output space $Y$ is taken to be the bounding space, here denoted $B$, of a FF of normed linear spaces $\{B_{s,t}\}$, or occasionally the corresponding extended space $B^e$. In general, the families $\{B_{s,t}\}$ and $\{A_{s,t}\}$ need not be the same. The notations for output spaces are analogous to those for input spaces. The comments about equivalence classes are valid for the $y_{s,t} \in B_{s,t}$.

A mapping $F : U \to Y$ is called a *global input-output mapping* (or usually just input-output mapping).

**Definition 4** *Let $(Y, F, U)$ be an input-output system. $F$ is a causal mapping and $(Y, F, U)$ is a causal system if and only if for all $t$ and for all $u, v \in U$ such that $\|u - v\|_t = 0$ it follows that $\|F(u) - F(v)\|_t = 0$.*

If $F$ satisfies this definition it determines a mapping from $U_t$ into $B_t$, denoted $\widetilde{F}_t$, that satisfies $\left\|\widetilde{F}_t u_t - (Fu)_t\right\|_t = 0$. We call $\widetilde{F}_t$ a *truncated input-output mapping* and define the *centered truncated input-output mapping* $F_t : U_0 \to Y_0$ by $F_t(u_0) \triangleq \mathrm{L}_t \widetilde{F}_t \mathrm{R}_t(u_0)$, where $\mathrm{R}_t \triangleq \mathrm{L}_{-t}$ is the right-shift by $t$. We assume that all systems in this report are causal. If $F$ satisfies this definition, it is causal in the usual sense. However, memory can affect causality.

**Definition 5** *Consider a left-expanded FF $\{A_{s,t}\}; -\infty \leq s < t < \infty$. The family $\{A_{s,t}\}$ and the norms $\|\cdot\|_t$ are said to be finite memory with memory length $M$ if there exists $0 < M < \infty$ such that $\|f\|_t = \|f\|_{t-M,t}$ for all $f \in L_0$, $t \in \Re$.*

Is is seen that if $w$ in (2-1) has finite support, that is $w(t) = 0$ for $t \geq M$, the weighted $L_p$ normed linear spaces are examples of finite memory spaces with memory length $M$. Finite memory may cause a system that is causal in the usual sense to not be causal according to Definition 4.

As in [13], the norms we use here for input-output mappings $F$, $\widetilde{F}_t$, and $F_t$ and for the natural states are the $N$-power norms, denoted $[\![ \cdot ]\!]_{(N)}$. (We omit the subscript $(N)$ when possible). Let $\Phi$ be a mapping from a normed linear space $X$ into a normed linear space $Z$. For any nonnegative integer $N$, the $N$-power norm for $\Phi$ is given by

$$[\![\Phi]\!]_{(N)} \triangleq \sup_{x \in X} \frac{\|\Phi(x)\|}{1 + \|x\|^N} \tag{2-5}$$

when the right side exists. We say $\Phi$ is *bounded (in N-power norm)* if $[\![\Phi]\!]_{(N)} < \infty$. If $\Phi$ is bounded, it carries bounded sets into bounded sets by the inequality

$$\|\Phi(x)\| \leq [\![\Phi]\!]_{(N)} \cdot \left(1 + \|x\|^N\right) .$$

However, boundedness of $\Phi$ does not in general imply continuity nor vice versa. The space of bounded operators from $U$ to $Y$ is denoted by $\mathcal{F}_N(U, Y)$. The space of bounded and continuous operators from $U$ to $Y$ is denoted by $\mathcal{C}_N(U, Y)$. We have chosen to use $N$-power norms rather than Lipschitz norms because they are less restrictive and because they are not so much influenced by the fine structure of a mapping. This last property seems to be important when one is dealing with an approximate system representation (see Heitman, [3] page 785). Other properties of these norms are given in the appendix of [13], which also gives comparisons with the Lipschitz norm. Although [13] uses the standard FFs, its appendix also applies to FFs in general. The $N$-power norms are special cases of a more general class of weighted supremum norms where $\|x\|^N$ in the denominator of the defining expression is replaced by an arbitrary, continuous, positive function of $x$; $X$ need not be normed (see Appendix A in [14]).

Using (2-5) on the truncated system mapping and the system mapping, we have

$$\|F_t\|_{(N)} = \sup_{u_0} \frac{\|F_t(u_0)\|_0}{1+\|u_0\|_0^N} = \sup_{u_t} \frac{\left\|\widetilde{F}_t(u_t)\right\|_t}{1+\|u_t\|_t^N} = \|\widetilde{F}_t\|_{(N)} \tag{2-6}$$

and

$$\|F\|_{(N)} = \sup_{u} \frac{\|F(u)\|}{1+\|u\|^N} \ . \tag{2-7}$$

Hence, we have the useful relationship

$$\|F\|_{(N)} \le \sup_t \|F_t\|_{(N)} = \sup_t \|\widetilde{F}_t\|_{(N)} \ . \tag{2-8}$$

Throughout this report, whenever there is reference to a system $(Y, F, U)$, the following three hypotheses are in effect unless specifically noted otherwise:

(A) The input space $U$ is either the bounding space $A$ of a FF of normed linear spaces that permits splicing or a shift-invariant subset of such an $A$. If $U$ is a proper subset of $A$, we require that it contain 0, but also that $u \in U$ implies both $u^t \mapsto 0_{t,\infty}$ and $0^t \mapsto u_{t,\infty}$ belong to $U$.

(B) The output space $Y$ is the bounding space $B$ of a FF of normed linear spaces.

(C) The global system operator $F$ satisfies Definition 4 (causality) with respect to the given $A$ and $B$.

The noncentered natural states $\widetilde{\xi}_t^u$ and the natural states $\xi_t^u$ are to be defined as operators with domains $\widetilde{\mathcal{D}}_t^u$ and $\mathcal{D}_t^u$, respectively, where

$$\widetilde{\mathcal{D}}_t^u \triangleq \{v_{t,\infty} \in U_{t,\infty} | u^t \mapsto v_{t,\infty} \in U, \ \forall u^t \subset u_t\}$$

and

$$\mathcal{D}_t^u \triangleq \{v_{0,\infty} \in U_{0,\infty} | u^t \mapsto \mathrm{R}_t v_{0,\infty} \in U, \ \forall u^t \subset u_t\} \ .$$

**Lemma 6** *Under the conditions just specified, given any $t \in \Re$ and any $u \in U$, there exists a mapping $\widetilde{\xi}_t^u : \widetilde{\mathcal{D}}_t^u \to Y_{t,\infty}$ such that $\left\|\widetilde{\xi}_t^u(v_{t,\infty}) - (F(u^t \mapsto v_{t,\infty}))_{t,\infty}\right\|_{t,\infty} = 0$ for all $v_{t,\infty} \in \widetilde{\mathcal{D}}_t^u$. Furthermore, if $u$ and $u'$ satisfy $\|u - u'\|_t = 0$ then $\widetilde{\xi}_t^u = \widetilde{\xi}_t^{u'}$. (Note that $\mathcal{D}_t^u = \mathcal{D}_t^{u'}$.)*

The proof of Lemma 6 is similar to the proof of Lemma 1 in [13]. Denote the Nerode equivalence class of defining inputs for $\widetilde{\xi}_t^u$ by $[u_t]_\xi$; i.e., $u'_t \in [u_t]_\xi$ if $\widetilde{\xi}_t^{u'} = \widetilde{\xi}_t^u$. Lemma 6 allows one to frame the following definition.

**Definition 7 ([10], [13], [14])** *The natural state for a system $(Y, F, U)$ induced at time $t$ by input $u_t$ is defined to be the operator $\xi_t^u$ from $\mathcal{D}_t^u$ to $Y_{0,\infty}$ given by*

$$\xi_t^u(v_{0,\infty}) = \mathrm{L}_t \widetilde{\xi}_t^u(\mathrm{R}_t v_{0,\infty}) \ . \tag{2-9}$$

*The set of natural states is denoted $\Xi$; the set of natural states that can be achieved at time $t$ is denoted $\Xi(t)$ so that $\Xi = \cup_{t\in\Re}\Xi(t)$. ($\Xi$ is referred to at times as the natural state space or the natural state set.)*

Obviously, the state set $\Xi$ is minimal since $\xi \in \Xi$ is defined as a mapping. Natural state may be defined for $F$ not necessarily causal; however, we only consider causal $F$ here.

The next hypothesis will often be needed but will not be in effect unless stated explicitly.

(D) The operators $F_t : U_0 \to Y_0$ are uniformly bounded in $N$-power norm for some fixed positive integer $N$ by a constant $C < \infty$ for all $t \in \Re$, and are an equicontinuous family of uniformly continuous mappings.

Hypothesis (D) gives that the global system operator is bounded with bound $C$ and uniformly continuous, see Lemmas 4 and 5 of [13].[3] Let $\mathcal{F}_N(X, Z)$ be the normed linear space of all mappings $\Phi : X \to Z$ with $\|\Phi\|_N < \infty$, and $\mathcal{C}_N(X, Z)$ be the normed linear subspace of $\mathcal{F}_N(X, Z)$ of all continuous $\Phi$. With Hypothesis (D) in force, $\Xi \subset \mathcal{C}_N(U_{0,\infty}, Y_{0,\infty})$; see Appendix A. Other properties that result from Hypothesis (D) are given in [13] and [14] and are mentioned at the time they are invoked. As stated previously, properties of the natural state and its relationship to $F$ are given in [10], [13], and [14].

[3] The proofs given in [13] were intended for a standard FF; however, such proofs usually hold for FFs in general if Condition 5 of Definition 1 is not used.

# 3 SOME STATE SPACE PROPERTIES

There are a few properties of the natural state space that can be related to constraints, or lack of constraints, on memory and on the class of admissible inputs. We have seen in [13] that the system determines the natural state space. We will find that some common systems, those represented by linear finite dimensional time-invariant differential equations, are determined by natural state space. That is, for these systems, the map from system to natural state space is one to one. Hence, the property that natural state space determines the system may be desirable for understanding system models in general. We examine how memory affects this property. The idea is that, in a sense, natural state defines the future part of the system. This, along with a tapered memory (a form of approximate finite memory, defined below), may specify the system. We give sufficient conditions, including tapered input spaces, by which systems are determined by natural state space. We also provide sufficient conditions for systems represented by polynomial integral operators to have this same property. We consider the effect memory has on reachability of states in the natural state space. It is shown that systems with finite memory input spaces or tapered input spaces exhibit a degree of reachability. Most of the results in this chapter are for time-invariant systems. We start with an example, which is very simple even though it is pathological, that points out certain negative possibilities.

**Example 8** *Let the basic function space $\mathcal{L}(\Re, \Re)$ for inputs be the set of real-valued functions on $\Re$ that are bounded and piecewise continuous (continuous except possibly at a finite number of points on any finite interval and with finite left- and right-hand limits at these points). The family of norms is given by $\|f\|_{s,t} = \sup_{s<\tau\le t}|f(\tau)|$, $f \in \mathcal{L}(\Re, \Re)$. The input space $U$ is either the bounding space $A$ for this family or a shift-invariant subset thereof. Let the basic function space for outputs be the Lebesgue space $L_\infty$ of real-valued functions. We want to have a finite observation interval of length $b > 0$ for the outputs, so the family of norms is given by*

$$\|f\|_{s,t} = \operatorname*{esssup}_{\max(s,t-b)<\tau\le t} |f(\tau)|, \ f \in L_\infty \ .$$

*The output space $Y$ is the bounding space for this family, and one has $Y_t = Y_{t-b,t}$. Let $T$, $T \in (0, \infty)$, be fixed and $h$ be a real-valued integrable function such that $h(s) = 0$ for $s \notin [0, T]$, and $\int_0^T |h(\tau)|\, d\tau \neq 0$.*

*Let $\overline{u} \triangleq \overline{\lim}_{t\to-\infty} u(t)$. Define a system mapping $F$ by*

$$[F(u)](s) \triangleq \int_{-\infty}^{s} h(s-\tau)u(\tau)d\tau + \overline{u}, \ -\infty < s < \infty \ . \tag{3-1}$$

*The output is $y_t \in Y_t$, where $y_t(s) = [\widetilde{F}_t(u_t)](s)$, $t-b < s \le t$, is given by (3-1) with $s$ restricted as indicated. The mapping $F$ is linear and time-invariant. Put $\mathbf{I} = \int_0^T |h(\tau)|d\tau$. Consider the continuity of $F_0$ and $F$:*

$$\left\|F_0(u) - F_0(u')\right\|_0 = \operatorname*{esssup}_{-b<s\le 0} \left| \int_{s-T}^{s} h(s-\tau)(u(\tau) - u'(\tau))d\tau + (\overline{u} - \overline{u}') \right|$$

$$\leq \mathbf{I}\cdot(\left\|u-u'\right\|_0)+\left|\overline{u}-\overline{u}'\right| \; . \tag{3-2}$$

*Since* $\|u-u'\|_0 < \epsilon$ *implies* $|\overline{u}-\overline{u}'| \leq \epsilon$, $F_0$ *is uniformly continuous and* $[\![F_0]\!]_N < \infty$ *for any* $N \geq 1$. *Also,* $F$ *is uniformly continuous and has finite N-power norm for any* $N \geq 1$.

*The state* $\xi_0^u$ *is given by*

$$[\xi_0^u(v_{0,\infty})](s) = \int_{\min(s-T,0)}^{0} h(s-\tau)u(\tau)d\tau + \overline{u} + \int_{\max(0,s-T)}^{s} h(s-\tau)v(\tau)d\tau, \;\; s \geq 0 \; , \tag{3-3}$$

*where it is to be noted that the first integral vanishes if* $s > T$. *Thus for* $s > T$, $[\xi_0^u(v_{0,\infty})](s)$ *depends on* $u$ *only through* $\overline{u}$. *Note that* $[\![\xi_0^u]\!]_N < \infty$ *for any* $N \geq 1$ *and is a uniformly continuous operator.*

*We consider this basic example now from two different aspects. The concept of reachability is needed in what follows.*

**Definition 9** *In a time-invariant system, the natural state* $\xi''$ *is reachable from* $\xi'$ *if there is an input that will drive the system from the state* $\xi'$ *to* $\xi''$ *in finite time.*

*We continue the example.*

*(a) To start with, it is immaterial whether* $U = A$ *or whether* $U$ *is bounded. Denote by* $U^1$ *and* $U^2$ *the sets of inputs* $u \in U$ *such that* $\overline{u} = 1$ *and* $\overline{u} = 2$ *respectively, and suppose they are not empty. Let* $\Xi^1$ *be the set of states* $\{\xi_0^u\}$, $u \in U^1$, *and* $\Xi^2$ *be the set of states* $\{\xi_0^w\}$, $w \in U^2$. *First, one may observe that* $\Xi^1$ *and* $\Xi^2$ *are disjoint sets, and indeed are actually a positive distance apart. In fact, for any* $u \in U^1$ *and* $w \in U^2$,

$$[\![\xi_0^u - \xi_0^w]\!]_N = \sup_{v_{0,\infty}} \frac{\|[\xi_0^u - \xi_0^w](v_{0,\infty})\|_{0,\infty}}{1+\|v\|_{0,\infty}^N}$$

$$\geq \frac{\|[\xi_0^u - \xi_0^w](v_{0,\infty})\|_{T,\infty}}{1+\|v\|_{0,\infty}^N} \tag{3-4}$$

*for any particular choice of* $v_{0,\infty}$. *Let* $v(s) = 0$ *for all* $s$. *Then, by (3-3), for this particular* $v$,

$$[\![\xi_0^u - \xi_0^w]\!]_N \geq \|[\xi_0^u - \xi_0^w](v_{0,\infty})\|_{T,\infty} \geq 1 \; .$$

*Furthermore, no state in* $\Xi^2$ *is reachable from any state in* $\Xi^1$ *and vice versa. Since the system is time-invariant, if the state* $\xi_0^u$, $u \in U^1$, *were to be reachable from the state represented by* $\xi_0^w$, $w \in U^2$, *there would have to be* $s < 0$ *and* $x \in U$ *such that* $\xi_0^{(R_s w)^s \mapsto x_{s,0}} = \xi_0^u$. *But this cannot be, because* $\overline{(R_s w)^s \mapsto x_{s,0}} = \overline{w}$, *no matter what* $s$ *and* $x$ *may be, implies* $(R_s w)^s \mapsto x_{s,0}$ *generates a state in* $\Xi^2$ *whereas* $\xi_0^u \in \Xi^1$. *Since* $\Xi^1$ *and* $\Xi^2$ *are a positive distance apart, we cannot even approximate* $\xi_0^u$ *starting from* $\xi_0^w$. *We note that it is admissible to consider a system as above with input space* $u \in U^1 \cup U^2$, *since this* $U$ *is translation-invariant.*

*(b) Consider now two systems. Two systems* $(Y, F^1, U)$ *and* $(Y, F^2, U)$ *(having the same input and output spaces) are identical if for all* $u \in U$ *we have* $y^1 = F^1(u) = y^2 = F^2(u)$. *Let* $\Xi$ *denote the natural state set for* $(Y, F^1, U)$, *and let* $H$ *denote the natural state set for* $(Y, F^2, U)$. *These two systems have identical natural state sets if, for all* $\xi \in \Xi$, *there exists* $\eta \in H$ *such that* *for all* $v_{0,\infty} \in U_{0,\infty}$ *we have* $y^1_{0,\infty} = \xi(v_{0,\infty}) = y^2_{0,\infty} = \eta\ (v_{0,\infty})$ *and vice versa. We set their input space to be* $U = A$*. We denote the first system by* $F^1$ *and define it by (3-1). The system mapping for the second system is given by*

$$[F^2(w)](s) = \int_{s-T}^{s} h(s-\tau)w(\tau)d\tau + 2 \cdot \overline{w}, \tag{3-5}$$

$-\infty < s < \infty$, $w \in U$*. These are clearly different systems. We claim, however, they have the same set of natural states. Let the state given by* $u$ *at* $t = 0$ *in the first system be* $\xi^u_0$ *and the state given by* $w$ *at* $t = 0$ *in the second system be* $\eta^w_0$*. For arbitrary* $u$*, we find* $w$ *such that*

$$[\xi^u_0(v_{0,\infty})](s) = [\eta^w_0(v_{0,\infty})](s) \tag{3-6}$$

*for all* $v$ *and all* $s \geq 0$*. From the formula for* $[\xi^u_0(v_{0,\infty})](s)$*, (3-3) and the corresponding formula for* $[\eta^w_0(v_{0,\infty})](s)$*, we see that any* $w \in U$*, such that (i)* $\overline{w} = 0.5 \cdot \overline{u}$ *and (ii)* $w(s) = u(s)$ *for* $-T < s \leq 0$ *will satisfy (3-6). The same sort of argument applies the other way around, so the two systems have the same set of natural states.*

Example 8 indicates some system behavior that can occur if there is "too much memory," meaning the present output depends on input from too far back in the past. Presumably one does not expect this sort of condition should occur, even in the model of a natural system. It does not happen with systems represented by the following types of linear differential equations.

**Example 10** *Consider two different systems* $\dot{x} = Ax + Bu$ *and* $\dot{x} = Fx + Gu$*, with output equation* $y = x$ *for both, where* $A \neq F$ *and where* $x(t) \in \Re^n$*,* $u(t) \in \Re^m$*. Let* $(A, B)$ *and* $(F, G)$ *be completely controllable. Then, the state vector* $x(t)$ *for each system may be any point in* $\Re^n$*. Let the natural state set for* $\dot{x} = Ax + Bu$ *be* $\Xi$ *and for* $\dot{x} = Fx + Gu$ *be* $H$*. From Example 2 of [13], a vector* $x(t) \in \Re^n$ *corresponds to a natural state* $\xi_t \in \Xi$ *and a natural state* $\eta_t \in H$*. However,* $\xi_t \neq \eta_t$*; e.g., consider the state transition property for* $\xi_t$ *and* $\eta_t$*. It may also be seen that* $\Xi \cap H = \emptyset$*. If* $\Xi = H$ *however, then* $A = F$ *and* $B = G$*. This also points out the difference between the state vector and the natural state.*

Using systems with tapered input spaces is a technical device which has been mentioned that, in some measure, at least prevents there being too much memory. It is used as a hypothesis in what follows. Part (a) of Example 8 illustrates how reachability and even approximate reachability can be prevented in a system with too much memory. We have no strong results on reachability, but there are some facts that can be deduced immediately and are summarized in the following remark and in Proposition 15. We will need the following definition.

**Definition 11 ([14])** *Let* $\mathcal{L}(\Re, E)$ *be a translation-invariant space of time functions from* $\Re$ *into* $E$ *with the property that, for each* $f \in \mathcal{L}$*, there is a real number* $T = T(f)$ *such that* $f(t) = 0$ *for all* $t < T$*. Then* $\mathcal{L}$ *is called a shifted-zero space.*

Shifted-zero input spaces may be convenient when considering unstable systems that do not "blow up" in a finite time or man-made systems that are activated at some finite time. For time-invariant systems, shifted-zero input spaces are equivalent to the system being defined over the half line. The following lemma gives sufficient conditions for the natural state space $\Xi$ to be connected. Connectedness is a weaker property than reachability for a natural state space because reachability implies connectedness.

**Lemma 12** *(Serakos, [14], Proposition 4.2.1.) Consider a (possibly time-varying) input-output system (Y,F,U). Assume that the state trajectories are continuous (e.g., use Proposition 4 of [13] or Proposition 4.1.2 of [14]) and that U is a shifted-zero space. Then $\Xi$ is arcwise connected.*

**Proof:** In fact, given $\xi'$ and $\xi''$, suppose $\xi' = \xi_t^u$ and $\xi'' = \xi_s^v$. Since U is a shifted-zero space, there is some $r < s, t$ such that $\|u - v\|_r = 0$. The state trajectory starting with $\xi_r^u$ generated by $u^r \!\mapsto v_{r,s}$ connects $\xi_r^u$ and $\xi''$. It then follows there is an arc running from $\xi'$ to $\xi_r^u$ and thence to $\xi''$. The states are arbitrary, so it follows that $\Xi$ is arcwise connected. (Note we are not saying that $\xi''$ is reachable from $\xi'$.) □

**Remark 13** *In a time-invariant system with a finite memory input space, it is obvious that any state is reachable from any other state, and the time required to reach a specified state is no more than the duration of the memory.*

From Remark 13, one would expect that a system that has an approximately finite memory input space — e.g., a system with a tapered input space — would exhibit approximate reachability, at least in the time-invariant case. But first, we present the definition of a tapered input space.

**Definition 14 ([10])** *Consider a left-expanded FF $\{A_{s,t}\}; -\infty \leq s < t < \infty$. For any $c > 0$ and $t$, let $G(c,t) \triangleq \{f \in \mathcal{L}_0, \|f\|_\tau \leq c$ for all $\tau \leq t\}$. Then $\{A_{s,t}\}$ and the norms $\|\cdot\|_t$ are said to be tapered if, for any $\epsilon > 0$, $c > 0$ and $t \in \Re$, there is $\delta = \delta(\epsilon, c, t) > 0$ such that $\|f\|_t \leq \|f\|_{t-\delta,t} + \epsilon$ for all $f \in G(c,t)$.*

Note that $\delta$ does not in fact depend on $t$ and also that $\|f\|_s \leq \|f\|_{s-\delta,s} + \epsilon$ for any $s \leq t$. When $w$ in (2-1) is integrable, the weighted $L_p$ normed linear spaces are examples of tapered FF of normed linear spaces. A specific class of tapered spaces are finite memory spaces. For the next proposition, we use the fact that Hypothesis D implies that the natural state is a uniformly continuous function of prior input (see Proposition 3 of [13]).

**Proposition 15** *(Professor Root provided this proposition.) Let Hypothesis D hold. Let $U$ be the bounding space and the associated norms be tapered. Then, given $\epsilon > 0$ and any two natural states $\xi'$ and $\xi''$, the natural state may be driven from $\xi'$ to approximate $\xi''$ as closely as desired. In fact, if $\xi' = \xi_0^x$ and $\xi'' = \xi_0^u$, there is $t < 0$ such that*

$$[\![\xi_0^{(\mathrm{R}_t x)^t \mapsto u_{t,0}} - \xi'']\!]_N \leq \epsilon \ .$$

**Proof:** Since $\xi' = \xi_0^x = \xi_t^{(\mathrm{R}_t x)}$ for any $t$, the result follows immediately from the definition of a tapered family, the fact $U$ is the bounding space, and since (see [13] and [14]) the natural state is a uniformly continuous function of prior input. □

A standard feature of the natural state space is that it is determined by the system. It is of mathematical interest to see when the map from system to natural state space is one to one.

Example 8 shows that natural state space is insufficient to determine the system. Example 10 gives a common example when it does. Hence, it may also be of practical interest (for instance in modeling man-made or natural systems) to resolve this same issue. Situations where natural state space determines the system are given by the following propositions and corollaries.

**Proposition 16** *([14], Proposition 4.2.3.) Consider two causal systems $(Y, F^1, U)$ and $(Y, F^2, U)$. Let these systems be time-invariant, continuous, and have the same state set $\Xi$. Then, if $U$ is bounded and tapered, they are identical systems.*

**Proof:** To show that $F^1$ and $F^2$ are identical, it is necessary to prove that $F^1(u) = F^2(u)$ for all $u \in U$. Denote a state for $F^1$ by $\xi$ and a state for $F^2$ by $\eta$. Given an input $u \in U$ to $F^2$, consider the state trajectory $t \to \eta_t^u$. For $T$ fixed, there exists $v \in U$, such that $\xi_T^v = \eta_T^u$ (for a time-invariant system, a natural state may be achieved at any time). Then,

$$\|\xi_T^v(x_{0,\infty}) - \eta_T^u(x_{0,\infty})\|_{0,\infty} = 0$$

for all $x_{0,\infty} \in \text{Domain}(\xi_T^v) = \text{Domain}(\eta_T^u)$, in particular,

$$\|\xi_T^v(\mathrm{L}_T(u_{T,\infty})) - \eta_T^u(\mathrm{L}_T(u_{T,\infty}))\|_{0,\infty}$$

$$= \left\|F^1(v^T \mapsto u_{T,\infty}) - F^2(u^T \mapsto u_{T,\infty})\right\|_{T,\infty} = 0 \ .$$

For $t = T-1$, there exists $w \in U$ such that $\xi_{T-1}^w = \eta_{T-1}^u$. Again

$$\left\|F^1(w^{T-1} \mapsto u_{T-1,\infty}) - F^2(u^{T-1} \mapsto u_{T-1,\infty})\right\|_{T-1,\infty} = 0 \ .$$

Therefore, it may be assumed that $v = w^{T-1} \mapsto u_{T-1,\infty}$. It now may be seen that a sequence of inputs $v^i \in U$, $i = 1, 2, \cdots$ may be constructed such that $\left\|v^i - u\right\|_{T+1-i,\infty} = 0$ with $\left\|F^1(v^i) - F^2(u)\right\|_{T+1-i,\infty} = 0$. Here, $v^1 = v^T \mapsto u_{T,\infty}$ and $v^2 = w^{T-1} \mapsto u_{T-1,\infty}$.

Since $U$ is bounded and tapered, given $\epsilon > 0$, there exists $\delta = \delta(\epsilon, T)$ such that

$$\left\|v^i - u\right\|_T \leq \left\|v^i - u\right\|_{T-\delta,T} + \epsilon$$

for all $i = 1, 2, \cdots$ (Definition 14). Take $i > \delta + 1$, then

$$\left\|v^i - u\right\|_T \leq 0 + \epsilon$$

since for such $i$,

$$\left\|v^i - u\right\|_{T-\delta,T} \leq \left\|v^i - u\right\|_{T+1-i,T} = 0 \ .$$

Thus $\overline{\lim}_{i\to\infty} \left\|v^i - u\right\|_T \leq \epsilon$. But $\epsilon$ is arbitrary, so $\lim_{i\to\infty} \left\|v^i - u\right\|_T = 0$, i.e., $v_T^i \to u_T$. Then,

$$\left\|F^1(u) - F^2(u)\right\|_T = \lim_{j\to\infty} \left\|F^1(u) - F^2(u)\right\|_{T-j,T}$$

$$= \lim_{j\to\infty} \left\|F^1(\lim_{i\to\infty} v^i) - F^2(u)\right\|_{T-j,T} = \lim_{j\to\infty} \lim_{i\to\infty} \left\|F^1(v^i) - F^2(u)\right\|_{T-j,T} \ .$$

But for fixed $j$, for all $i \geq j+1$

$$\left\|F^1(v^i) - F^2(u)\right\|_{T-j,T} = 0 \ .$$

Hence, for all $j$,

$$\lim_{i\to\infty} \left\|F^1(v^i) - F^2(u)\right\|_{T-j,T} = 0$$

and thus,

$$\lim_{j\to\infty} \lim_{i\to\infty} \left\|F^1(\nu^i) - F^2(u)\right\|_{T-j,T} = 0 \ . \qquad \square$$

In the following two corollaries, additional conditions are used to apply Proposition 16 to systems with unbounded input spaces and time-varying systems.

**Corollary 17** *For an input-output system $(Y, F, U)$ with state space $\Xi$ and unbounded input space $U$, let*

$$\Xi^C = \{\xi_t^u \in \Xi | \ \|u\| \leq C\} \ .$$

*Let $(Y, F^1, U)$ have state space $\Xi$ and $(Y, F^2, U)$ have state space $H$, and as in Proposition 16, be time-invariant continuous systems with $U$ tapered but unbounded. Then the two systems are identical if $\Xi^C = H^C$ for all $C < \infty$.*

**Proof:** Apply Proposition 16 to all designated $u \in U$ with $C = \|u\|$. $\square$

A statement may be made concerning Proposition 16 and time-varying systems if one degree of freedom is removed in a manner similar to the previous corollary.

**Corollary 18** *For an input-output system $(Y, F, U)$ with state space $\Xi$, let*

$$\Xi(t) = \{\xi_t^u \in \Xi\} \ .$$

*Let $(Y, F^1, U)$ having state space $\Xi$ and $(Y, F^2, U)$ having state space $H$ be continuous but time-varying systems. Let the input space $U$ be bounded and tapered. Then the two systems are identical if $\Xi(t) = H(t)$ for all $t \in \Re$.*

**Proof:** Trace through the proof of Proposition 16 applying the new hypothesis where appropriate. $\square$

Next, we look at the question, does state determine system for systems represented by polynomial integral operators? These are an important class of systems because they may be used to approximate classes of systems; see [8]. A scalar version of the following proposition is in [16]. We provide a definition of a polynomial integral operator.

**Definition 19 ([14])** *By a time-invariant $N^{th}$ degree causal, polynomial integral operator with $M$-dimensional input space $u = [u_1, \cdots, u_M] \in U$ with norm $\|u\|_{s,t} = \max_i \|u_i\|_{s,t}$ and $\|u_i\|_{s,t} = \sup_{s<\tau\leq t} |u_i(\tau)|$ and $P$-dimensional output space with norm $\|y\|_{s,t} = \max_i \|y_i\|_{s,t}$ and $\|y_i\|_{s,t} = \sup_{\max(s,t-b)<\tau\leq t} |y_i(\tau)|$, $b \geq 0$, we mean a mapping of the form:*

$$[F_p(u)](t) = \sum_{n=1}^{N} [F_p^n(u)](t), \ \ p \in [1, 2, \cdots, P] \tag{3-7}$$

*where*

$$[F_p^n(u)](t) = \int_0^\infty \cdots \int_0^\infty \sum_{i_1,\cdots,i_n=1}^{M} f_p^{i_1,\cdots,i_n}(\sigma_1,\cdots,\sigma_n)u_{i_1}(t-\sigma_1)\cdots u_{i_n}(t-\sigma_n)d\sigma_1,\cdots,d\sigma_n \ .$$

*The above sum is taken over all combinations without repeating; hence, there are $M^n$ terms. Such an operator is unchanged if the kernels $f_p^{i_1,\cdots,i_n}$, all $p$ and $n$ are symmetrized. The symmetrized kernel of $f_p^{i_1,\cdots,i_n}$, denoted by $\widetilde{f}_p^{i_1,\cdots,i_n}$, is defined by*

$$\widetilde{f}_p^{i_1,\cdots,i_n}(\sigma_1,\cdots,\sigma_n) \triangleq \frac{1}{n!}\sum_{\pi} f_p^{i_{\pi(1)},\cdots,i_{\pi(n)}}(\sigma_{\pi(1)},\cdots,\sigma_{\pi(n)}) \ .$$

*(The sum is taken over all permutations $\pi$, total $n!$ terms).*

Advantages and disadvantages of polynomial integral operators are discussed in Section 3 of Root, [8]. It turns out that their generality is cited for both of these considerations. Their generality makes them useful while, at the same time, can make them unwieldy for some applications.

**Proposition 20** *Consider an input-output system represented by a polynomial integral operator. The following conditions are imposed on the kernels $f_p^{i_1,\cdots,i_n}$:*
*(a) $f_p^{i_1,\cdots,i_n}(\tau_1,\cdots,\tau_n)$ is continuous on $\Re^n$.*
*(b) $f_p^{i_1,\cdots,i_n}(\tau_1,\cdots,\tau_n)$ is absolutely integrable on $\Re^n$.*
*(c) The symmetrized kernel $\widetilde{f}$ satisfies*

$$\sup_{\tau_1,\cdots,\tau_i} \widetilde{f}_p^{i_1,\cdots,i_n}(\tau_1,\cdots,\tau_i,\tau_{i+1},\cdots,\tau_n)$$

*is absolutely integrable as a function of $\tau_{i+1},\cdots,\tau_n$ for all $1 \le i < n$, for all $1 \le n \le N$ and all $p \in P$.*

*Denote the class of operators given by (3-7) with kernels $f_p^{i_1,\cdots,i_n}$ satisfying (a), (b), and (c) by $V(N)$ and denote the subclass of $V(N)$ with symmetric kernels by $V_s(N)$. Under these conditions,*

*(i) In the class of systems $(Y, F, U)$ with $U$ and $Y$ as given immediately above and with $F \in V_s(N)$, the natural state set $\Xi$ determines the kernels $\widetilde{f}_1,\cdots,\widetilde{f}_N$.*

*(ii) In the class of systems $(Y, F, U)$ with $U$ and $Y$ as in (i) and $F \in V(N)$, the natural state set $\Xi$ determines the system, i.e., determines $F$.*

**Proof:** Given in Appendix B.

Consider systems $(Y, F, U)$ with $F \in V(N)$. Since each $u \in U$ is bounded, $F(u) \in Y$ by condition (b); hence, these systems are well-defined. The assertion (ii) follows directly from (i). In fact, if two systems of the specified class have the same state set $\Xi$, they are represented by the same set of symmetrized kernels; hence, the (perhaps unsymmetrized) kernels used in their description produce the same mapping. Our proof of assertion (i) goes by first taking an interval around a point where the two kernels are assumed to be different. By looking at the form of the natural state, an input is selected to show that in fact the two kernels have to agree at that point. Initially, scalar $N$-degree kernels are considered. Then $N$-order scalar polynomial kernels are considered. Finally, vector polynomial kernels are considered.

We note that we need only conditions (a) and (b) if the systems are homogeneous integral operators, but we need (a), (b), and (c) for the general polynomial case. Incidentally, (i) implies that under the stated conditions for any $F \in V(N)$, the symmetrized kernels are determined uniquely — a fact we presume known.

The justification for this proposition is that the restriction of a bounded input space used in Proposition 16 or the restrictions on the state spaces used in Corollaries 17 or 18 are not required for Proposition 20. However, the conclusion here is weaker in that the uniqueness of the system mapping is shown to hold only within the special class of systems defined by polynomial integral operators.

# 4 STATE DIFFERENTIABILITY PROPERTIES IN INPUT-OUTPUT SYSTEMS

In this chapter, differentiability properties of various maps involving the natural state of an input-output system are analyzed. The analysis is patterned after the continuity analysis in Section 3 of [13] with the objective of obtaining differentiability analogs to some of the continuity results. In this regard, these results may be viewed as upgrades of results presented in [13]. These results are inherited from appropriate differentiability assumptions on the input-output system. Gateaux, as well as Frechet, differentials and derivatives are used in these results as they appear in [14]. Here we present the Frechet case only. Time-varying systems are considered in this chapter.

We present three results. All require the differentiability of the input-output system. The first gives the differentiability of the natural state. The second gives the differentiability of the map from past input to natural state. This result requires the Frechet derivative remainder converges uniformly to zero. The last result implies the differentiability of the natural state trajectories. This result requires that $U$ be shift differentiable and that $t \to F_t$ be a shift differentiable system trajectory.

This chapter may be contrasted with the work of Jakubczyk [4], [5], and also with Sussmann, [17]. In comparing [4], [5], and [17] with this report, we must note the distinction between state vectors and natural states. Example 2 of [13] and Example 10 (this report) aid in noting the distinction. In [4], [5], and [17] state vectors are considered, while in this report natural states are considered. In [4] and [5], a finite dimensional state space differential equation realization of an input-output system is developed. This is in contrast to the differentiability results for the natural states developed here, which includes differentiability of the natural state trajectories. The assumptions used in [4] and [5] are necessary and sufficient for providing a finite dimensional state variable realization of an input-output system. A fundamental feature of these papers is that the state vectors are points in an $n$-dimensional differentiable manifold. The differential equation representation derived in these papers is of minimum dimension for a given input-output system. That is, the state vectors are of minimum dimension (or $n$ is minimized). In [17] a differential equation representation of minimal dimension is determined for a system that is already in state variable or differential form but not necessarily minimal. The natural state space in this report is a subset of an infinite dimensional Banach space (e.g., it is not necessarily linear, locally Euclidean, or a manifold.) As mentioned in the introduction, the natural state space has the essential state properties. The natural state space is minimal as set. This is different than minimizing the dimension of the state vector. Also, a state vector in [4] and [5] is contained in an open set; however, a natural state may not be contained in an open set (although Remark 13 gives conditions when the natural state space is connected). No coordinate chart can be defined on the natural state space since its dimension is not specified. Dimensionality (finite or not) is not relevant in this report. Differentiability of individual natural state space trajectories is considered. Because of these reasons, [4], [5], and [17] are not comparable with this report.

In this chapter we will assume that the domains of the states are large enough for the stated operations to be defined (e.g., $U = A$, see Hypothesis A).

**Definition 21** *The input-output system $(Y, F, U)$ has a Frechet derivative if for $u, v \in U$, with $u$ contained in an open subset of $U$,*

$$F(u+v) = F(u) + L(u,v) + W(u,v) \tag{4-1}$$

*where $L(u,\cdot) : U \to Y$ is linear and bounded and $W(u,\cdot) : U \to Y$ has the property that*

$$\lim_{\|v\|\to 0} \frac{\|W(u,v)\|}{\|v\|} = 0 \ .$$

*$L(u,v)$ is called the Frechet differential of $(Y, F, U)$ at $u$ in the direction $v$, and $W(u,v)$ is called the remainder. The operator $L(u,\cdot) : U \to Y$ is called the Frechet derivative of $F$ at $u$ and is denoted $F'(u)$, [2].*

We make a special note that Frechet differentiability implies continuity, see Cheney, [2]. The following proposition is a differentiability upgrade to Proposition 1 in [13].

**Proposition 22** *If the Frechet derivative of $(Y, F, U)$ at $u$ exists, then the Frechet derivative of the natural state map $\xi_t^u : U_{0,\infty} \to Y_{0,\infty}$ at $\mathrm{L}_t u_{t,\infty}$ exists.*

Proof: Let $v_{0,\infty} = \mathrm{L}_t u_{t,\infty}$, $w_{0,\infty} \in U_{0,\infty}$. We have

$$\begin{aligned}
&\xi_t^u(v_{0,\infty} + w_{0,\infty}) - \xi_t^u(v_{0,\infty}) \\
&= [\mathrm{L}_t F(u^t \!\longleftrightarrow \mathrm{R}_t(v_{0,\infty} + w_{0,\infty}))]_{0,\infty} - [\mathrm{L}_t F(u^t \!\longleftrightarrow \mathrm{R}_t v_{0,\infty})]_{0,\infty} \\
&= [\mathrm{L}_t F(u^t \!\longleftrightarrow \mathrm{R}_t v_{0,\infty} + 0^t \!\longleftrightarrow \mathrm{R}_t w_{0,\infty}) - \mathrm{L}_t F(u^t \!\longleftrightarrow \mathrm{R}_t v_{0,\infty})]_{0,\infty} \\
&= [\mathrm{L}_t L(u^t \!\longleftrightarrow \mathrm{R}_t v_{0,\infty}, 0^t \!\longleftrightarrow \mathrm{R}_t w_{0,\infty}) + \mathrm{L}_t W(u^t \!\longleftrightarrow \mathrm{R}_t v_{0,\infty}, 0^t \!\longleftrightarrow \mathrm{R}_t w_{0,\infty})]_{0,\infty} \\
&= [\mathrm{L}_t L(u^t \!\longleftrightarrow \mathrm{R}_t v_{0,\infty}, 0^t \!\longleftrightarrow \mathrm{R}_t w_{0,\infty})]_{0,\infty} + [\mathrm{L}_t W(u^t \!\longleftrightarrow \mathrm{R}_t v_{0,\infty}, 0^t \!\longleftrightarrow \mathrm{R}_t w_{0,\infty})]_{0,\infty} \ .
\end{aligned} \tag{4-2}$$

Obviously, $L^1(v_{0,\infty},\cdot) \stackrel{\Delta}{=} [\mathrm{L}_t L(u^t \!\longleftrightarrow \mathrm{R}_t v_{0,\infty}, 0^t \!\longleftrightarrow \mathrm{R}_t \cdot)]_{0,\infty} : U_{0,\infty} \to Y_{0,\infty}$ inherits linearity from $L(u,\cdot)$. Also,

$$\begin{aligned}
[]L^1[] &= \sup_{w_{0,\infty}} \frac{\left\|L^1(v_{0,\infty}, w_{0,\infty})\right\|_{0,\infty}}{\|w_{0,\infty}\|_{0,\infty}} = \sup_{w_{0,\infty}} \frac{\left\|\mathrm{L}_t L(u^t \!\longleftrightarrow \mathrm{R}_t v_{0,\infty}, 0^t \!\longleftrightarrow \mathrm{R}_t w_{0,\infty})\right\|_{0,\infty}}{\|w_{0,\infty}\|_{0,\infty}} \\
&= \sup_{w_{0,\infty}} \frac{\left\|\mathrm{L}_t L(u^t \!\longleftrightarrow \mathrm{R}_t v_{0,\infty}, 0^t \!\longleftrightarrow \mathrm{R}_t w_{0,\infty})\right\|_{0,\infty}}{\|0^t \!\longleftrightarrow \mathrm{R}_t w_{0,\infty})\|} \cdot \frac{\left\|0^t \!\longleftrightarrow \mathrm{R}_t w_{0,\infty})\right\|}{\|w_{0,\infty}\|_{0,\infty}} \\
&\leq \sup_{w_{0,\infty}} \frac{\left\|\mathrm{L}_t L(u^t \!\longleftrightarrow \mathrm{R}_t v_{0,\infty}, 0^t \!\longleftrightarrow \mathrm{R}_t w_{0,\infty})\right\|}{\|0^t \!\longleftrightarrow \mathrm{R}_t w_{0,\infty})\|} \leq []L[] \ ;
\end{aligned}$$

hence, $L^1$ is bounded. Defining $W^1(v_{0,\infty}, w_{0,\infty}) \stackrel{\Delta}{=} [\mathrm{L}_t W(u^t \!\longleftrightarrow \mathrm{R}_t v_{0,\infty}, 0^t \!\longleftrightarrow \mathrm{R}_t w_{0,\infty})]_{0,\infty}$, we have

$$\lim_{\|w_{0,\infty}\|_{0,\infty}\to 0} \frac{\left\|W^1(v_{0,\infty}, w_{0,\infty})\right\|_{0,\infty}}{\|w_{0,\infty}\|_{0,\infty}} = \lim_{\|w_{0,\infty}\|_{0,\infty}\to 0} \frac{\left\|\mathrm{L}_t W(u^t \!\longleftrightarrow \mathrm{R}_t v_{0,\infty}, 0^t \!\longleftrightarrow \mathrm{R}_t w_{0,\infty})\right\|_{0,\infty}}{\|w_{0,\infty}\|_{0,\infty}}$$

$$\leq \lim_{\|0^t \leftrightarrow \mathrm{R}_t w_{0,\infty}\| \to 0} \frac{\left\| W(u^t \leftrightarrow \mathrm{R}_t v_{0,\infty}, 0^t \leftrightarrow \mathrm{R}_t w_{0,\infty}) \right\|}{\|0^t \leftrightarrow \mathrm{R}_t w_{0,\infty}\|} = 0 \ .$$

This last inequality comes from the change in the time interval of the norm in the numerator. Hence, the Frechet derivative of $\xi_t^u$ at $\mathrm{L}_t u_{t,\infty}$ exists with Frechet derivative $L^1$ and remainder $W^1$. As a final point, it appears that $L^1$ as it is defined above depends on $u^t$. We shall show that any element in the Nerode equivalence class of $\xi_t^u$, i.e., any $u_t \in [u_t]_\xi$ may be used. Let $u_t, u_t' \in [u_t]_\xi$. Then

$$[\![ L^1((v_{0,\infty}, \cdot) - L'^1(v_{0,\infty}, \cdot) ]\!] = \sup_{w_{0,\infty}} \frac{\left\| L^1(v_{0,\infty}, w_{0,\infty}) - L'^1(v_{0,\infty}, w_{0,\infty}) \right\|_{0,\infty}}{\|w_{0,\infty}\|_{0,\infty}}$$

$$= \sup_{w_{0,\infty}} \frac{\left\| [\mathrm{L}_t L(u^t \leftrightarrow \mathrm{R}_t v_{0,\infty}, 0^t \leftrightarrow \mathrm{R}_t w_{0,\infty})]_{0,\infty} - [\mathrm{L}_t L(u'^t \leftrightarrow \mathrm{R}_t v_{0,\infty}, 0^t \leftrightarrow \mathrm{R}_t w_{0,\infty})]_{0,\infty} \right\|_{0,\infty}}{\|w_{0,\infty}\|_{0,\infty}}$$

$$= \sup_{w_{0,\infty}} \lim_{c \to 0} \frac{\left\| \begin{array}{c} [\mathrm{L}_t L(u^t \leftrightarrow \mathrm{R}_t v_{0,\infty}, 0^t \leftrightarrow c \cdot \mathrm{R}_t w_{0,\infty})]_{0,\infty} \\ -[\mathrm{L}_t L(u'^t \leftrightarrow \mathrm{R}_t v_{0,\infty}, 0^t \leftrightarrow c \cdot \mathrm{R}_t w_{0,\infty})]_{0,\infty} \end{array} \right\|_{0,\infty}}{c \cdot \|w_{0,\infty}\|_{0,\infty}}$$

$$= \sup_{w_{0,\infty}} \lim_{c \to 0} \frac{\left\| \begin{array}{c} [\mathrm{L}_t W(u'^t \leftrightarrow \mathrm{R}_t v_{0,\infty}, 0^t \leftrightarrow c \cdot \mathrm{R}_t w_{0,\infty}) \\ -\mathrm{L}_t W(u^t \leftrightarrow \mathrm{R}_t v_{0,\infty}, 0^t \leftrightarrow c \cdot \mathrm{R}_t w_{0,\infty})]_{0,\infty} \end{array} \right\|_{0,\infty}}{c \cdot \|w_{0,\infty}\|_{0,\infty}} = 0 \ . \qquad \square$$

The following proposition is a differentiability analog to Proposition 3 in [13].

**Proposition 23** *Let all the natural states have the same domain. If the Frechet derivative of a bounded map $F$ at $u$ exists, if $\alpha_U = +\infty$ and if $(\|W(u,v)\| / \|v\|) \to 0$ as $\|v\| \to 0$ uniformly in $u$, then the Frechet derivative of the map $S : U_t \to \Xi$ defined by $S(v_t) = \xi_t^v$ at $u_t$ exists.*

Proof: Calculate

$$\begin{aligned} &\xi_t^{u+v}(w_{0,\infty}) - \xi_t^u(w_{0,\infty}) \\ &= [\mathrm{L}_t F((u+v)^t \leftrightarrow \mathrm{R}_t w_{0,\infty}) - \mathrm{L}_t F(u^t \leftrightarrow \mathrm{R}_t w_{0,\infty})]_{0,\infty} \\ &= [\mathrm{L}_t F(u^t \leftrightarrow \mathrm{R}_t w_{0,\infty} + v^t \leftrightarrow 0_{t,\infty}) - \mathrm{L}_t F(u^t \leftrightarrow \mathrm{R}_t w_{0,\infty})]_{0,\infty} \\ &= [\mathrm{L}_t L(u^t \leftrightarrow \mathrm{R}_t w_{0,\infty}, v^t \leftrightarrow 0_{t,\infty}) + \mathrm{L}_t W(u^t \leftrightarrow \mathrm{R}_t w_{0,\infty}, v^t \leftrightarrow 0_{t,\infty})]_{0,\infty} \ . \end{aligned}$$

First, consider the map $u_t \to \xi_t^u(w_{0,\infty})$ for $w_{0,\infty}$ fixed. Clearly, $[\mathrm{L}_t L(u^t \leftrightarrow \mathrm{R}_t w_{0,\infty}, \cdot \leftrightarrow 0_{t,\infty})]_{0,\infty} : U_0 \to Y_{0,\infty}$ is linear, and

$$\lim_{\|v_t\|_t \to 0} \frac{\left\| \mathrm{L}_t W(u^t \leftrightarrow \mathrm{R}_t w_{0,\infty}, v^t \leftrightarrow 0_{t,\infty}) \right\|_{0,\infty}}{\|v_t\|_t}$$

$$\leq \lim_{\|v^t \leftrightarrow 0_{t,\infty}\|_t \to 0} \frac{\left\| W(u^t \leftrightarrow \mathrm{R}_t w_{0,\infty}, v^t \leftrightarrow 0_{t,\infty}) \right\|}{K_U^{-1} \cdot \|v^t \leftrightarrow 0_{t,\infty}\|} = 0$$

by increasing the time interval of the norm in the numerator and by the assumption on $U$. (Note that since $\alpha_U = +\infty$, $\|v_t\|_t \to 0$ iff $\left\| v^t \leftrightarrow 0_{t,\infty} \right\| \to 0$.) Now, consider $S$. Define an operator $\Lambda(u_t,\cdot) : U_t \to \mathcal{C}_N(U_{0,\infty}, Y_{0,\infty})$ for some positive integer $N$ by

$$[\Lambda(u_t, v_t)](w_{0,\infty}) \triangleq [\mathrm{L}_t L(u^t \hookrightarrow \mathrm{R}_t w_{0,\infty}, v^t \hookrightarrow 0_{t,\infty})]_{0,\infty} \ .$$

Note that $\Lambda$ is linear in $v_t$. Also, define an operator $\Omega(u_t, v_t) : U_t \times U_t \to \mathcal{C}_N^\infty(U_{0,\infty}, Y_{0,\infty})$ by

$$[\Omega(u_t, v_t)](w_{0,\infty}) \triangleq [\mathrm{L}_t W(u^t \hookrightarrow \mathrm{R}_t w_{0,\infty}, v^t \hookrightarrow 0_{t,\infty})]_{0,\infty} \ .$$

It follows from the uniformity hypothesis on $W$ that

$$\lim_{\|v_t\|_t \to 0} \frac{[\![\Omega(u_t, v_t)]\!]}{\|v_t\|_t} = \lim_{\|v_t\|_t \to 0} \frac{\sup\limits_{w_{0,\infty}} \dfrac{\left\|\mathrm{L}_t W(u^t \hookrightarrow \mathrm{R}_t w_{0,\infty}, v^t \hookrightarrow 0_{t,\infty})\right\|_{0,\infty}}{1 + \|w_{0,\infty}\|_{0,\infty}^N}}{\|v_t\|_t}$$

$$\leq \sup_{w_{0,\infty}} \lim_{\|v^t \hookrightarrow 0_{t,\infty}\| \to 0} \frac{\left\|W(u^t \hookrightarrow \mathrm{R}_t w_{0,\infty}, v^t \hookrightarrow 0_{t,\infty})\right\|_{t,\infty}}{K_U^{-1} \cdot \|v^t \hookrightarrow 0_{t,\infty}\|} = 0 \ .$$

If $S$ is bounded, then $\Lambda$ is bounded, (see Theorem 2, Chapter 3 of [2]).

$$[]S[] = \sup_{u_t} \frac{[]\xi_t^u[]}{1 + \|u_t\|_t^N} = \sup_{u_t} \frac{\sup\limits_{v} \dfrac{\|\xi_t^u(v_{0,\infty})\|_{0,\infty}}{1 + \|v_{0,\infty}\|_{0,\infty}^N}}{1 + \|u_t\|_t^N}$$

$$\leq \sup_{u,v} \frac{\left\|F(u^t \hookrightarrow \mathrm{R}_t v_{0,\infty})\right\|}{\left(1 + \|u_t\|_t^N\right)\left(1 + \|v_{0,\infty}\|_{0,\infty}^N\right)} = \sup_{u,v} \frac{\left\|F(u^t \hookrightarrow \mathrm{R}_t v_{0,\infty})\right\|}{1 + \|u_t\|_t^N + \|v_{0,\infty}\|_{0,\infty}^N + \|u_t\|_t^N \|v_{0,\infty}\|_{0,\infty}^N}$$

$$\leq (2^N + 1) \cdot \sup_{u,v} \frac{\left\|F\left(u^t \hookrightarrow \mathrm{R}_t v_{0,\infty}\right)\right\|}{1 + \left(\|u^t\|_t + \|v_{0,\infty}\|_{0,\infty}\right)^N} \leq (2^N + 1) \cdot \sup_{u,v} \frac{\left\|F\left(u^t \hookrightarrow \mathrm{R}_t v_{0,\infty}\right)\right\|}{1 + \|u^t \hookrightarrow \mathrm{R}_t v_{0,\infty}\|^N}$$

$$= (2^N + 1) \cdot \sup_u \frac{\|F(u)\|}{1 + \|u\|^N} = (2^N + 1) \cdot []F[] \ .$$

Hence, $S$ has Frechet derivative $\Lambda(u_t, \cdot)$ at $u_t$ with remainder $\Omega(u_t, \cdot)$. □

We will need the following two definitions in the next proposition, which is a differentiability upgrade (analog) to Proposition 4 in [13].

**Definition 24** *The input space $U$ is shift-differentiable if, for all $u \in U$,*

$$\mathrm{L}_h u = u + h \cdot d + e(h) \tag{4-3}$$

*where $d, e(h) \in U$ and $\lim_{h \to 0} \|e(h)\|_t / h = 0$ for all $t \in \Re$. $d$ is the shift derivative of $u$.*

A shift-differentiable input space $U$ is also shift-continuous, as mentioned in Appendix C. A shift-differentiable input $u$ does not necessarily have a derivative with respect to time (at least not everywhere; see Appendix C.)

**Definition 25** *A system trajectory is a shift-differentiable system trajectory $t \to F_t$ if*

$$\mathrm{L}_h F \mathrm{R}_h = F + h \cdot G + H(h) \tag{4-4}$$

*where $G \in \mathcal{C}_N^\infty(U, Y)$ and $H(h) \in \mathcal{C}_N^\infty(U, Y)$ and such that*

$$\lim_{h \to 0} \frac{[]H(h)_t[]}{h} = 0 \text{ for all } t \in \Re\cdot \tag{4-5}$$

Conditions that include a semigroup generator for the system trajectory for a shift-differentiable system trajectory are given in Appendix D. Appendix D also gives conditions so that (4-5) is uniform in $t$.

**Proposition 26** *Let all the natural states have the same domain. Let Hypothesis D from Page 2-6 hold. Assume that the input space $U$ and the system trajectory $t \to F_t$ are shift-differentiable and that the Frechet derivative $L(u, \cdot)$ of $F$ exists for all $u \in U$. In addition, with regard to the shift differentiability of $U$, assume that $\lim_{h \to 0} \|e(h)\|_s / h \to 0$ uniformly for $s \leq T$ for all $T \in \Re$ for each $u \in U$. Also, assume that (4-5) is uniform in $t$. Then for all $u \in U$, the state trajectory $t \to \xi_t^u$ is differentiable.*

Proof: First calculate:

$$\xi_{t+h}^u(v_{0,\infty}) - \xi_t^u(v_{0,\infty}) = [\mathrm{L}_{t+h} F(u^{t+h} \leftrightarrow \mathrm{R}_{t+h} v_{0,\infty}) - \mathrm{L}_t F(u^t \leftrightarrow \mathrm{R}_t v_{0,\infty})]_{0,\infty}$$

$$= [\mathrm{L}_{t+h} F(u^{t+h} \leftrightarrow \mathrm{R}_{t+h} v_{0,\infty}) - \mathrm{L}_t F((\mathrm{L}_h u)^t \leftrightarrow \mathrm{R}_t v_{0,\infty})]_{0,\infty}$$
$$+[\mathrm{L}_t F(u^t \leftrightarrow \mathrm{R}_t v_{0,\infty} + (hd + e(h))^t \leftrightarrow 0_{t,\infty}) - \mathrm{L}_t F(u^t \leftrightarrow \mathrm{R}_t v_{0,\infty})]_{0,\infty} \; . \tag{4-6}$$

So then,

$$\lim_{h \to 0} \frac{\xi_{t+h}^u(v_{0,\infty}) - \xi_t^u(v_{0,\infty})}{h}$$
$$= \lim_{h \to 0} \frac{[\mathrm{L}_{t+h} F(u^{t+h} \leftrightarrow \mathrm{R}_{t+h} v_{0,\infty}) - \mathrm{L}_t F((\mathrm{L}_h u)^t \leftrightarrow \mathrm{R}_t v_{0,\infty})]_{0,\infty}}{h}$$
$$+ \lim_{h \to 0} \frac{\begin{array}{c}[\mathrm{L}_t(L(u^t \leftrightarrow \mathrm{R}_t v_{0,\infty}, (hd + e(h))^t \leftrightarrow 0_{t,\infty}) \\ +W(u^t \leftrightarrow \mathrm{R}_t v_{0,\infty}, (hd + e(h))^t \leftrightarrow 0_{t,\infty}))]_{0,\infty}\end{array}}{h} \; . \tag{4-7}$$

Consider part of the second term on the right-hand side of (4-7):

$$\left\| \lim_{h \to 0} \frac{[\mathrm{L}_t(W(u^t \leftrightarrow \mathrm{R}_t v_{0,\infty}, (hd + e(h))^t \leftrightarrow 0_{t,\infty}))]_{0,\infty}}{h} \right\|_{0,\infty}$$

$$\leq \lim_{h \to 0} \frac{\left\| W(u^t \leftrightarrow \mathrm{R}_t v_{0,\infty}, (hd + e(h))^t \leftrightarrow 0_{t,\infty}) \right\|}{\|(hd + e(h))^t \leftrightarrow 0_{t,\infty}\|} \cdot \frac{\left\|(hd + e(h))^t \leftrightarrow 0_{t,\infty}\right\|}{h} = \mathbf{I} \; . \tag{4-8}$$

Now,

$$\lim_{h \to 0} \frac{\left\|(hd + e(h))^t \leftrightarrow 0_{t,\infty}\right\|}{h} = \lim_{h \to 0} \sup_s \frac{\left\|(hd + e(h))^t \leftrightarrow 0_{t,\infty}\right\|_s}{h} = \mathbf{II} \; .$$

To compute this limit, there are two cases to consider:

(i) $s > t$:

$$\mathbf{II} \le \|d\| + \lim_{h\to 0} \frac{\|e(h)\|_t}{h} = \|d\| \ .$$

(ii) $s \le t$:

$$\mathbf{II} \le \lim_{h\to 0} \sup_{s<t} \left( \frac{\|hd\|_s + \|e(h)\|_s}{h} \right) \le \|d\| + \lim_{h\to 0} \sup_{s<t} \frac{\|e(h)\|_s}{h} = \|d\| \ .$$

Therefore, **I** is zero.

Consider the first term on the right-hand side of (4-7):

$$\lim_{h\to 0} \frac{[\mathrm{L}_{t+h} F(u^{t+h} \!\!\longleftrightarrow \mathrm{R}_{t+h} v_{0,\infty}) - \mathrm{L}_t F((\mathrm{L}_h u)^t \!\!\longleftrightarrow \mathrm{R}_t v_{0,\infty})]_{0,\infty}}{h}$$

$$= \lim_{h\to 0} \frac{[\mathrm{L}_t(\mathrm{L}_h F(\mathrm{R}_h((\mathrm{L}_h u)^t \!\!\longleftrightarrow \mathrm{R}_t v_{0,\infty})) - F((\mathrm{L}_h u)^t \!\!\longleftrightarrow \mathrm{R}_t v_{0,\infty}))]_{0,\infty}}{h}$$

$$= \lim_{h\to 0} \frac{[\mathrm{L}_t(\mathrm{L}_h F(\mathrm{R}_h((u+hd+e(h))^t \!\!\longleftrightarrow \mathrm{R}_t v_{0,\infty})) - F((u+hd+e(h))^t \!\!\longleftrightarrow \mathrm{R}_t v_{0,\infty}))]_{0,\infty}}{h}$$

$$= \lim_{h\to 0} \frac{[\mathrm{L}_t(h \cdot G((u+hd+e(h))^t \!\!\longleftrightarrow \mathrm{R}_t v_{0,\infty}) + H(h)((u+hd+e(h))^t \!\!\longleftrightarrow \mathrm{R}_t v_{0,\infty}))]_{0,\infty}}{h} \ .$$

where $G$ and $H$ are from Definition 25. Now,

$$\lim_{h\to 0} \frac{[\mathrm{L}_t(h \cdot G((u+hd+e(h))^t \!\!\longleftrightarrow \mathrm{R}_t v_{0,\infty}))]_{0,\infty}}{h} = [\mathrm{L}_t G(u^t \!\!\longleftrightarrow \mathrm{R}_t v_{0,\infty})]_{0,\infty} \ . \qquad (4\text{-}9)$$

And

$$\lim_{h\to 0} \frac{\left\|[\mathrm{L}_t(H(h)((u+hd+e(h))^t \!\!\longleftrightarrow \mathrm{R}_t v_{0,\infty}))]_{0,\infty}\right\|_{0,\infty}}{h}$$

$$\le \lim_{h\to 0} \sup_{\tau} \frac{\left\|H(h)((u+hd+e(h))^t \!\!\longleftrightarrow \mathrm{R}_t v_{0,\infty})\right\|_\tau}{h}$$

$$\le \lim_{h\to 0} \sup_{\tau} \frac{[]H(h)_\tau[]}{h} \cdot (1 + \left\|(u+hd+e(h))^t \!\!\longleftrightarrow \mathrm{R}_t v_{0,\infty}\right\|^N) = 0 \ .$$

Hence, from (4-7), (4-8), and (4-9), the derivative of the state trajectory $t \to \xi^u_t$ is

$$\lim_{h\to 0} \frac{\xi^u_{t+h}(\cdot) - \xi^u_t(\cdot)}{h} = [\mathrm{L}_t G(u^t \!\!\longleftrightarrow \mathrm{R}_t \cdot)]_{0,\infty} + [\mathrm{L}_t(L(u^t \!\!\longleftrightarrow \mathrm{R}_t \cdot, d^t \!\!\longleftrightarrow 0_{t,\infty}))]_{0,\infty} \ .\square \qquad (4\text{-}10)$$

Thus, for time-invariant systems, differentiability of the natural state trajectory follows from the general conditions and the shift differentiability of the input FF. Note that $u^t$ appears in (4-10). However, similar to Proposition 22, any element in the Nerode equivalence class of $u^t$ may be used. Indeed, consider a time $T$ and natural state $\xi^u_T$ and suppose there is $u' \in U$ such that $\xi_T^{u'} = \xi_T^u$. Then for all $t \geq T$, $v \in U$ by the state transition property,

$$\xi_t^{u^T \mapsto v_{T,\infty}} = \xi_t^{u'^T \mapsto v_{T,\infty}} \ .$$

If the natural state trajectories are differentiable, then the derivatives of the natural state trajectories are identical for all $t \geq T$. Hence,

$$\frac{d(\xi_t^{u^T \mapsto v_{T,\infty}})}{dt} = \frac{d(\xi_t^{u'^T \mapsto v_{T,\infty}})}{dt} \ ,$$

for all $t \geq T$, where for $t = T$, the derivative is taken from the right. Assume that these two inputs are shift differentiable, then for all $t \geq T$

$$[\mathrm{L}_t G(u'^T \mapsto v^{T,t} \mapsto \mathrm{R}_t \cdot)]_{0,\infty} + [\mathrm{L}_t (L(u'^T \mapsto v^{T,t} \mapsto \mathrm{R}_t \cdot, d'^t \mapsto 0_{t,\infty})]_{0,\infty}$$
$$= [\mathrm{L}_t G(u^T \mapsto v^{T,t} \mapsto \mathrm{R}_t \cdot)]_{0,\infty} + [\mathrm{L}_t (L(u^T \mapsto v^{T,t} \mapsto \mathrm{R}_t \cdot, d^t \mapsto 0_{t,\infty})]_{0,\infty} \ .$$

For $t = T$, limits are taken from the right in computing these terms. Hence, it does not matter if $u', d'$ or $u, d$ are used in (4-10).

From Proposition 26 we get a differential equation representation for the natural state trajectory:

$$\frac{d\xi_t^u}{dt} = [\mathrm{L}_t G(u^t \mapsto \mathrm{R}_t \cdot)]_{0,\infty} + [\mathrm{L}_t (L(u^t \mapsto \mathrm{R}_t \cdot, d^t \mapsto 0_{t,\infty}))]_{0,\infty} \ . \tag{4-11}$$

Note that the domains and ranges of the two terms on the right in (4-11) have the same form as the domain and range of a natural state.

**Example 27** *Consider[4] the causal input-output system described by a second degree polynomial integral operator of the form*

$$y(r) = [F(u)](r) = \int_0^\infty \int_0^\infty f(r, \tau_1, \tau_2) u(r - \tau_1) u(r - \tau_2) d\tau_1 d\tau_2, \ -\infty < r < \infty \tag{4-12}$$

*where, purely for simplicity, the input time functions, the kernel $f$, and hence the output time functions are scalar valued (otherwise the product forming the integrand must be defined as a tensor product). Obviously such a system is time-varying unless the kernel $f$ is not actually a function of its first argument. It will be assumed that $f(r, \tau_2, \tau_2)$ is symmetric in $\tau_1$ and $\tau_2$. This entails no loss of generality since symmetrizing $f$ does not change the value of the integral.[5] We choose ordinary $L_2$ norms for the FF of the input space and the $L_\infty$ norms for*

[4] *This example is a continuation of Example 1, [13].*

[5] *From the definition of a symmetrized kernal, $f_{sym}(r, \tau_1, \tau_2)) = \frac{1}{2}(f(r, \tau_1, \tau_2) + f(r, \tau_2, \tau_1))$. Hence, the symmetry property, $f_{sym}(r, \tau_1, \tau_2) = f_{sym}(r, \tau_2, \tau_1)$. To show that using the regular (unsymmetrized) or symmetrized kernal within the integral (4-12) is the same, we have $\int_{-\infty}^{\infty} \int_{-\infty}^{\infty} f_{sym}(r, \tau_1, \tau_2) d\tau_1 d\tau_2 = \int_{-\infty}^{\infty} \int_{-\infty}^{\infty} \frac{1}{2}(f(r, \tau_1, \tau_2) + f(r, \tau_2, \tau_1)) d\tau_1 d\tau_2 = \frac{1}{2} \int_{-\infty}^{\infty} \int_{-\infty}^{\infty} f(r, \tau_1, \tau_2) d\tau_1 d\tau_2 + \frac{1}{2} \int_{-\infty}^{\infty} \int_{-\infty}^{\infty} f(r, \tau_2, \tau_1) d\tau_2 d\tau_1 = \int_{-\infty}^{\infty} \int_{-\infty}^{\infty} f(r, \tau_1, \tau_2) d\tau_1 d\tau_2$ .*

*the FF of the output space. Thus, e.g.,*

$$\|u\|_0 = \left( \int_{-\infty}^{0} |u(t)|^2 \, dt \right)^{1/2} \tag{4-13}$$

*and*

$$\|y\|_0 = \operatorname*{ess\,sup}_{-\infty < t \leq 0} |y(t)| \ . \tag{4-14}$$

*Finally, assume that for all s there is $M < \infty$ such that*

$$\int_0^{\infty} \int_0^{\infty} |f(s, \tau_1, \tau_2)|^2 d\tau_1 d\tau_2 \leq M \ . \tag{4-15}$$

*We now determine the form of the operators $F_t$ and $\xi_t^u$. First, $[\widetilde{F}_t(u_t)](r) = [F(u)](r)$ for $r \leq t$, so $\widetilde{F}_t$ is given by Equation (4-12) with $-\infty < r \leq t$. Then,*

$$[F_t(u_0)](p) = [\mathrm{L}_t \widetilde{F}_t(\mathrm{R}_t u_0)](p) = [F(\mathrm{R}_t u)](t+p)$$

$$= \int_0^{\infty} \int_0^{\infty} f((t+p), \tau_1, \tau_2) u(p - \tau_1) u(p - \tau_2) d\tau_1 d\tau_2, \ -\infty < p \leq 0 \ . \tag{4-16}$$

*The natural state $\xi_t^u$ is given by*

$$[\xi_t^u(v_{0,\infty})](\sigma) = [\mathrm{L}_t \widetilde{\xi}_t^u(\mathrm{R}_t v_{0,\infty})](\sigma)$$

$$= [F(u_t \longmapsto R_t v_{0,\infty})](t+\sigma), \ \sigma \geq 0 \ . \tag{4-17}$$

*Evaluation of the right side of (4-17) yields*

$$[\xi_t^u(v_{0,\infty})](\sigma) = \int_0^{\infty} \int_0^{\infty} f(t+\sigma, \sigma+\mu_1, \sigma+\mu_2) u(t-\mu_1) u(t-\mu_2) d\mu_1 d\mu_2$$

$$+2 \int_0^{\sigma} \int_0^{\infty} f(t+\sigma, \sigma+\mu_1, \tau_2) u(t-\mu_1) v(\sigma-\tau_2) d\mu_1 d\tau_2$$

$$+ \int_0^{\sigma} \int_0^{\sigma} f(t+\sigma, \tau_1, \tau_2) v(\sigma-\tau_1) v(\sigma-\tau_2) d\tau_1 d\tau_2, \ \sigma \geq 0 \ . \tag{4-18}$$

*We need to consider the boundedness and continuity of the system operator. Looking at the issue of boundedness: From (4-12), (4-13), and (4-15) it follows for all $t$ and all $p \leq 0$ that*

$$|[F_t(u_0)](p)| = \left| \int_0^{\infty} \int_0^{\infty} f((t+p), \tau_1, \tau_2) u(p-\tau_1) u(p-\tau_2) d\tau_1 d\tau_2 \right|$$

$$\leq \int_0^{\infty} \int_0^{\infty} |f((t+p), \tau_1, \tau_2) u(p-\tau_1) u(p-\tau_2)| \, d\tau_1 d\tau_2 \ ,$$

*using the Schwarz inequality,*

$$\leq \left(\int_0^\infty \int_0^\infty |f((t+p),\tau_1,\tau_2)|^2 \, d\tau_1 d\tau_2\right)^{1/2} \left(\int_0^\infty \int_0^\infty |u(p-\tau_1)u(p-\tau_2)|^2 \, d\tau_1 d\tau_2\right)^{1/2}$$

$$\leq M \cdot \left(\int_0^\infty \int_0^\infty |u(p-\tau_1)|^2 \, d\tau_1\right)^{1/2} \left(\int_0^\infty \int_0^\infty |u(p-\tau_2)|^2 \, d\tau_2\right)^{1/2}$$

$$\leq M \cdot \|u\|_0^2 \tag{4-19}$$

*From (4-19),*

$$\|F_t\|_N = \sup_{u \in U_0} \frac{\|F_t(u_0)\|_0}{1+\|u\|_0^N} \leq M \sup_{u \in U_0} \frac{\|u\|_0^2}{1+\|u\|_0^N} \ ,$$

*so that if $N \geq 2$, the $F_t$ are all bounded by $M$. Then by (2-8), $\|F\| \leq M$. By Appendix A, or Proposition 2 of [13], for each $u \in U$, $\|\xi_t^u\|$ is bounded in $t$. $U_0$ may be any subset of $A_0$, including $A_0$ itself. Turning to the issue of continuity,*

$$\left|[F_t(u_0) - F_t(u_0')](p)\right|$$

$$= \left|\int_0^\infty \int_0^\infty f((t+p),\tau_1,\tau_2)\left(u(p-\tau_1)u(p-\tau_2) - u'(p-\tau_1)u'(p-\tau_2)\right) d\tau_1 d\tau_2\right|$$

$$\leq \int_0^\infty \int_0^\infty \Big| f((t+p),\tau_1,\tau_2)\left\{u(p-\tau_1)\left(u(p-\tau_2)-u'(p-\tau_2)\right)\right.$$

$$\left. + \left(u(p-\tau_1)-u'(p-\tau_1)\right)u'(p-\tau_2)\right\}\Big| \, d\tau_1 d\tau_2 \ ,$$

*using the Schwarz inequality,*

$$\leq \left(\int_0^\infty \int_0^\infty |f((t+p),\tau_1,\tau_2)|^2 \, d\tau_1 d\tau_2\right)^{1/2}$$

$$\cdot \left(\int_0^\infty \int_0^\infty \left|u(p-\tau_1)\left(u(p-\tau_2)-u'(p-\tau_2)\right) + \left(u(p-\tau_1)-u'(p-\tau_1)\right)u'(p-\tau_2)\right|^2 d\tau_1 d\tau_2\right)^{1/2}$$

*using the Minkowski inequality,*

$$\leq \left(\int_0^\infty \int_0^\infty |f((t+p),\tau_1,\tau_2)|^2 \, d\tau_1 d\tau_2\right)^{1/2}$$

$$\cdot \left(\left(\int_0^\infty \int_0^\infty \left|u(p-\tau_1)\left(u(p-\tau_2)-u'(p-\tau_2)\right)\right|^2 d\tau_1 d\tau_2\right)^{1/2}\right.$$

$$\left. + \left(\int_0^\infty \int_0^\infty \left|\left(u(p-\tau_1)-u'(p-\tau_1)\right)u'(p-\tau_2)\right|^2 d\tau_1 d\tau_2\right)^{1/2}\right)$$

$$\leq \left(\int_0^\infty \int_0^\infty |f((t+p),\tau_1,\tau_2)|^2 \, d\tau_1 d\tau_2\right)^{1/2}$$

$$\cdot\left(\left(\int_0^\infty |u(p-\tau_1)|^2\,d\tau_1\right)^{1/2}\left(\int_0^\infty \left|u(p-\tau_2)-u'(p-\tau_2)\right|^2 d\tau_2\right)^{1/2}\right.$$

$$\left.+\left(\int_0^\infty \left|u(p-\tau_1)-u'(p-\tau_1)\right| d\tau_1\right)^{1/2}\left(\int_0^\infty \left|u'(p-\tau_2)\right|^2 d\tau_2\right)^{1/2}\right)$$

$$\leq M\cdot\left\|u_0-u_0'\right\|_0\cdot\left(\|u_0\|_0+\left\|u_0'\right\|_0\right)\ ; \tag{4-20}$$

*therefore, $F_t$ is continuous. However, if $U_0$ is unbounded we cannot guarantee that $F_t$ is uniformly continuous. So now suppose that $U_0$ is bounded, $\|u_0\|_0 \leq M'$, $u_0 \in U_0$. Then by (4-20)*

$$\left\|F_t(u_0)-F_t(u_0')\right\|_0=\sup_{p\leq 0}\left|[F_t(u_0)-F_t(u_0')](p)\right|\leq M\cdot 2M'\left\|u_0-u_0'\right\|_0$$

*so that $\{F_t\}$ is an equicontinuous family of uniformly continuous operators. This gives that $F$ is uniformly continuous (shown in Lemma 4 of [13]). Then from Appendix A, the natural states $\xi_t^u$ are uniformly continuous. Furthermore, it follows that the conclusions of Proposition 3 of [13] (the mapping input to natural state is uniformly continuous) and Proposition 5 of [13] (future natural state depends continuously on past natural state) hold. From Proposition 4 of [13] the natural state trajectory is continuous if the system trajectory $t \to F_t$ is continuous (e.g., the system is time invariant).*

*We next consider the differentiability properties presented in this chapter. To start with, we consider whether the system (4-12) has a Frechet derivative. We have*

$$F\left(u+v\right)=\int_0^\infty\int_0^\infty f(r,\tau_1,\tau_2)\left(\left(u+v\right)\left(r-\tau_1\right)\right)\left(\left(u+v\right)\left(r-\tau_2\right)\right)d\tau_1 d\tau_2$$

$$=\int_0^\infty\int_0^\infty f(r,\tau_1,\tau_2)u(r-\tau_1)u(r-\tau_2)d\tau_1 d\tau_2$$

$$+2\cdot\int_0^\infty\int_0^\infty f(r,\tau_1,\tau_2)u(r-\tau_1)v(r-\tau_2)d\tau_1 d\tau_2$$

$$+\int_0^\infty\int_0^\infty f(r,\tau_1,\tau_2)v(r-\tau_1)v(r-\tau_2)d\tau_1 d\tau_2\ . \tag{4-21}$$

*We compare the terms of (4-21) with those in Definition 21. Looking at (4-12),*

$$F\left(u\right)=\int_0^\infty\int_0^\infty f(r,\tau_1,\tau_2)u(r-\tau_1)u(r-\tau_2)d\tau_1 d\tau_2\ .$$

*Define*

$$L\left(u,v\right)\triangleq 2\cdot\int_0^\infty\int_0^\infty f(r,\tau_1,\tau_2)u(r-\tau_1)v(r-\tau_2)d\tau_1 d\tau_2 \tag{4-22}$$

*and we see $L\left(u,\cdot\right)$ is linear. $L$ is bounded as is seen from*

$$[]L\left(u,\cdot\right)[]=\sup_v\frac{\left\|L\left(u,v\right)\right\|}{\|v\|}$$

$$= \sup_v \frac{\sup_r \left|2 \cdot \int_0^\infty \int_0^\infty f(r,\tau_1,\tau_2) u(r-\tau_1) v(r-\tau_2) d\tau_1 d\tau_2\right|}{\|v\|}$$

$$\leq \sup_v \frac{2 \cdot \sup_r \sqrt{\int_0^\infty \int_0^\infty \|f(r,\tau_1,\tau_2)\|^2 d\tau_1 d\tau_2} \sqrt{\int_0^\infty \int_0^\infty \|u(\tau_1) v(\tau_2)\|^2 d\tau_1 d\tau_2}}{\|v\|}$$

$$= \sup_v \frac{2 \cdot \sup_r \sqrt{\int_0^\infty \int_0^\infty \|f(r,\tau_1,\tau_2)\|^2 d\tau_1 d\tau_2} \sqrt{\int_0^\infty \|u(\tau_1)\|^2 d\tau_1} \sqrt{\int_0^\infty \|v(\tau_2)\|^2 d\tau_2}}{\|v\|}$$

$$\leq \frac{2M \cdot \|u\| \, \|v\|}{\|v\|} \leq 2M \cdot \|u\| \ .$$

*Define*

$$W(u,v) \triangleq \int_0^\infty \int_0^\infty f(r,\tau_1,\tau_2) v(r-\tau_1) v(r-\tau_2) d\tau_1 d\tau_2 \ . \tag{4-23}$$

*Then*

$$\lim_{\|v\| \to 0} \frac{\|W(u,v)\|}{\|v\|} = \lim_{\|v\| \to 0} \frac{\sup_r \left|\int_0^\infty \int_0^\infty |f(r,\tau_1,\tau_2) v(r-\tau_1) v(r-\tau_2)| \, d\tau_1 d\tau_2\right|}{\|v\|}$$

$$\leq \lim_{\|v\| \to 0} \frac{M \cdot \|v\|^2}{\|v\|} = 0 \ .$$

*Hence the input-output system has a Frechet derivative and this decides the applicability of Proposition 22. We see that (4-23) is uniform in u, and therefore, Proposition 23 is applicable. Looking at Proposition 26, a number of conditions need to be satisfied for this proposition to give differentiable state trajectories. The shift differentiability of $L_2$ is considered in Appendix C. The sufficient conditions of Proposition C-5 provide that $u \in L_2$ is shift differentiable if it has a time derivative and has a bounded remainder, as in (C-6). Our system trajectory is shift differentiable if it is time invariant (this is a sufficient condition, other specific cases may be considered.) These conditions are sufficient for Proposition 26 to apply to this input-output system. Checking into the form of (4-11), if the system is time-invariant, G is zero. Then*

$$\frac{d\xi_t^u}{dt} = [\mathrm{L}_t(L(u^t \mapsto \mathrm{R}_t \cdot, d^t \mapsto 0_{t,\infty}))]_{0,\infty} \ ,$$

*where d is the shift derivative of u. Substituting (4-22),*

$$\frac{d\xi_t^u}{dt} (v_{0,\infty})(r)$$

$$= 2 \cdot \int_0^\infty \int_0^\infty f(t+r,\tau_1,\tau_2) \left(u^t \mapsto \mathrm{R}_t v_{0,\infty}\right)(t+r-\tau_1) \left(d^t \mapsto 0_{t,\infty}\right)(t+r-\tau_2) d\tau_1 d\tau_2 \ ,$$

*where $0 \leq r < \infty$.*

# 5 SUMMARY AND CONCLUSIONS

It is pointed out that the natural state set is unique to an input-output system if a tapered input space is used (Proposition 16) or if the system has a polynomial integral operator representation (Proposition 20). (A form of this result may also hold for the differential equation representation in Chapter 4.) The uniqueness of the natural state space to an input-output system may be applicable to the modeling of physical systems and may be of mathematical interest. That is, if uniqueness of the natural state space to a system is evident, then steps should be taken to incorporate this feature in a mathematical model. This is a property of common system representations and hence may already be tacitly assumed. For instance, linear finite dimensional time-invariant differential systems are uniquely determined by their natural state space (Example 10). Example 8 shows that systems modeled with excessive memory may not uniquely define a natural state space. This same example along with Proposition 15 showed that reachability is affected by memory length. Sufficient conditions for a connected natural state space are given in Lemma 12. Connectedness and reachability of the natural state space may also be properties of interest in modeling. We gave conditions in [13] for the map from system to natural state space to be well-defined.

It is also pointed out that various differentiability properties of the natural state and natural state trajectory are inherited from the input-output system. These results are differentiability upgrades to the continuity results presented in [13]. Proposition 22 gives that if the Frechet derivative of the input-output system at $u$ exists, then the Frechet derivative of the natural state at $\mathrm{L}_t u_{t,\infty}$ exists. Proposition 23 gives conditions for the Frechet derivative of the map from input up to time $t$ to natural state at time $t$ to exist. Proposition 26 gives conditions such that the Frechet derivative exists along a natural state trajectory. This proposition has the result of a differential equation representation for the natural state trajectory (Equation (4-11)). These results, as well as their Gateaux derivative versions, mostly appear in [14]. The differentiability results presented in this report also may be considered as aids in modeling of physical systems.

# APPENDIX A

# BOUNDEDNESS AND CONTINUITY OF THE NATURAL STATE

In this appendix, we demonstrate that the natural state inherits boundedness and continuity from the input-output system. NSWCDD's publications format guide requires that appendices be self-contained, so we repeat here the definitions of fitted families (FFs) of time functions and other relevant definitions from the main text of this report.

**Definition A-1** *([A-1]) Let* $\mathcal{L} = \mathcal{L}(\Re, E)$ *be a linear space of time functions from* $\Re$ *into a Banach space* $E$ *such that any translate of a function in* $\mathcal{L}$ *is also a function in* $\mathcal{L}$*. Let* $\mathcal{N} = \{\|\cdot\|_{s,t}, -\infty < s < t < \infty\}$ *be a family of seminorms on* $\mathcal{L}$ *satisfying the following conditions:*

*(1) For* $f_1, f_2 \in \mathcal{L}$*, if* $f_1(\tau) = f_2(\tau)$ *for* $s < \tau \leq t$ *then* $\|f_1 - f_2\|_{s,t} = 0$*.*
*(2) Let* $\mathrm{L}_\tau$ *denote shift to the left by* $\tau$*. For all* $f \in \mathcal{L}$*,* $\|\mathrm{L}_\tau f\|_{s-\tau,t-\tau} = \|f\|_{s,t}$*.*
*(3) Let* $r < s < t$*. Then for all* $f \in \mathcal{L}$*,* $\|f\|_{s,t} \leq \|f\|_{r,t}$*.*
*(4) Let* $r < s < t$*. Then for all* $f \in \mathcal{L}$*,* $\|f\|_{r,t} \leq \|f\|_{r,s} + \|f\|_{s,t}$*.*
*(5) There exists* $0 < \alpha \leq \infty$ *and* $K \geq 1$ *such that if* $0 < t - r \leq \alpha$ *and* $r < s < t$*, then for all* $f \in \mathcal{L}$*,* $\|f\|_{r,s} \leq K \|f\|_{r,t}$*.*

*The pair* $(\mathcal{L}, \mathcal{N})$ *is called a FF of seminorms on* $\mathcal{L}$*. The normed linear space formed from equivalence classes of functions in* $\mathcal{L}$ *with norm* $\|\cdot\|_{s,t}$ *is denoted* $A_{s,t}$*. The elements of* $A_{s,t}$ *are the equivalence classes determined by:* $f \sim g$*,* $f, g \in \mathcal{L}$ *if and only if* $\|f - g\|_{s,t} = 0$*. They are denoted* $u_{s,t}$*,* $y_{s,t}$*, etc. The set* $\{A_{s,t}\}$*,* $-\infty < s < t < \infty$*, is the FF of normed linear spaces given by* $(\mathcal{L}, \mathcal{N})$*.*

For $f \in \mathcal{L}$, put $\|f\|^{s,t} \triangleq \sup_{s<\tau\leq t} \|f\|_{s,\tau}$. A FF $(\mathcal{L}, \mathcal{N})$ and $\{A_{s,t}\}$, $-\infty < s < t < \infty$, can be augmented to include $\|\cdot\|_{-\infty,t}$ by taking the limit $s \to -\infty$, since by (3) of Definition A-1 $\|f\|_{s,t}$ is monotone nondecreasing as $s \to -\infty$ with $t$ fixed. Let $\mathcal{L}_0 = \{f \in \mathcal{L} |\ \lim_{s\to-\infty} \|f\|_{s,t} < \infty, t \in \Re\}$. For $f \in \mathcal{L}_0$, define

$$\|f\|_t \triangleq \lim_{s\to-\infty} \|f\|_{s,t} = \|f\|_{-\infty,t} . \tag{A-1}$$

With the meaning of $(\mathcal{L}, \mathcal{N})$ thus extended, $\|\cdot\|_{s,t}$ is defined for $-\infty \leq s < t < \infty$. The *left-expanded FF of seminorms* is thereby defined and is denoted $(\mathcal{L}_0, \mathcal{N})$. It still satisfies all the Conditions $(1), \cdots, (5)$.

We next define $\|\cdot\|_{s,\infty}$ and $A_{s,\infty}$. For a FF, this is done by taking the supremum. Let $\mathcal{L}_{00} = \{f \in \mathcal{L}_0 |\ \sup_t \|f\|_t < \infty\}$. For $f \in \mathcal{L}_{00}$ define

$$\|f\|_{s,\infty} \triangleq \sup_{t>s} \|f\|_{s,t} \ ; \ -\infty \leq s . \tag{A-2}$$

It may be readily verified that if $(\mathcal{L}, \mathcal{N})$ is a FF for indices satisfying $-\infty < s < t < \infty$ then, with definitions given by (A-1) and (A-2), $(\mathcal{L}_{00}, \mathcal{N})$ is a FF for indices satisfying $-\infty \leq s < t < \infty$ and satisfies Conditions 1, 2, 3, and 5 of Definition A-1 for indices $-\infty \leq s < t \leq \infty$. $\{(\mathcal{L}_{00}, \mathcal{N}),\ \|\cdot\|_{s,t},\ -\infty \leq s < t \leq \infty\}$ is called the *expanded family of seminorms* determined by $(\mathcal{L}, \mathcal{N})$.

For $f \in \mathcal{L}_{00}$, we put

$$\|f\| \triangleq \sup_{t\in\Re} \|f\|_t = \|f\|_{-\infty,\infty} . \tag{A-3}$$

The normed linear space consisting of equivalence classes of functions in $\mathcal{L}_{00}$ with the norm (A-3) is called the *bounding space* $A$ for the family $\{A_{s,t}\}$. For $-\infty \leq r < s < t < \infty$, and $h$,

$g \in \mathcal{L}$, the *splice* of $h$ and $g$ over $(r, t]$ at $s$ is defined and equals $f$ if

$$f(\tau) = \begin{cases} h(\tau), & r < \tau \leq s \\ g(\tau), & s < \tau \leq t \end{cases}$$

belongs to $\mathcal{L}$. It is denoted $f_{r,t} = h_{r,s} \mapsto g_{s,t}$. For $t = \infty$, the splice of $h$ and $g$ equals $f$ if

$$f(\tau) = \begin{cases} h(\tau), & r < \tau \leq s \\ g(\tau), & s < \tau \end{cases}$$

belongs to $\mathcal{L}$. It is denoted $f_{r,\infty} = h^{r,s} \mapsto g_{s,\infty}$. Let $\Phi$ be a mapping from a normed linear space $X$ into a normed linear space $Z$. As shown in Appendix A of [A-2], for any nonnegative integer $N$, the $N$-power norm for $\Phi$ is given by

$$\|\Phi\|_{(N)} \stackrel{\Delta}{=} \sup_{x \in X} \frac{\|\Phi(x)\|}{1 + \|x\|^N} \tag{A-4}$$

when the right side exists.

Boundedness and continuity transfer from the input-output system to the natural state. Consider the boundedness of the natural state.

$$\begin{aligned} \|\xi_t^u\| &= \sup_{v_{0,\infty}} \frac{\|\xi_t^u (v_{0,\infty})\|_{0,\infty}}{1 + \|v_{0,\infty}\|_{0,\infty}^N} = \sup_{v_{t,\infty}} \frac{\left\|\widetilde{\xi}_t^u (v_{t,\infty})\right\|_{t,\infty}}{1 + \|v_{t,\infty}\|_{t,\infty}^N} \\ &= \sup_{v_{t,\infty}} \frac{\left\|F\left(u^t \mapsto v_{t,\infty}\right)\right\|_{t,\infty}}{1 + \|u^t \mapsto v_{t,\infty}\|^N} \cdot \frac{1 + \left\|u^t \mapsto v_{t,\infty}\right\|^N}{1 + \|v_{t,\infty}\|_{t,\infty}^N} \\ &\leq \sup_{v_{t,\infty}} \frac{\left\|F\left(u^t \mapsto v_{t,\infty}\right)\right\|}{1 + \|u^t \mapsto v_{t,\infty}\|^N} \cdot \sup_{v_{t,\infty}} \frac{1 + \left\|u^t \mapsto v_{t,\infty}\right\|^N}{1 + \|v_{t,\infty}\|_{t,\infty}^N} \\ &\leq \|F\| \cdot \sup_{v_{t,\infty}} \frac{1 + \sup_s \left\|u^t \mapsto v_{t,\infty}\right\|_s^N}{1 + \|v_{t,\infty}\|_{t,\infty}^N} = \mathbf{I} \;. \end{aligned}$$

There are two cases of values for the variable "$s$" in $\mathbf{I}$ to consider.

(i) Suppose $s \leq t$:

$$\mathbf{I} \leq \|F\| \cdot \left(1 + \|u\|^N\right)$$

(ii) Suppose $s > t$:

$$\mathbf{I} \leq \|F\| \cdot \sup_{v_{t,\infty}} \frac{1 + \sup_{s>t} \left\|u^t \mapsto v_{t,\infty}\right\|_s^N}{1 + \|v_{t,\infty}\|_{t,\infty}^N}$$

$$\leq [\![F]\!] \cdot \left( 1 + \sup_{v_{t,\infty}} \frac{\left( \|u_t\|_t + \|v_{t,\infty}\|_{t,\infty} \right)^N}{1 + \|v_{t,\infty}\|_{t,\infty}^N} \right)$$

$$= [\![F]\!] \cdot \left( 1 + \sup_{v_{t,\infty}} \frac{\|u\|^N + a_{N-1} \|u\|^{N-1} \|v_{t,\infty}\|_{t,\infty} + \cdots + a_0 \|v_{t,\infty}\|_{t,\infty}^N}{1 + \|v_{t,\infty}\|_{t,\infty}^N} \right)$$

where $a_0, \cdots, a_{N-1}$ are constants obtained from expanding the numerator of the previous equation. This last expression is less than or equal to

$$\leq [\![F]\!] \cdot \left( 1 + \|u\|^N + a_{N-1} \|u\|^{N-1} + \cdots + a_0 \right) \ .$$

Consider the continuity of the natural state.

$$\|\xi_t^u (v_{0,\infty}) - \xi_t^u (w_{0,\infty})\|_{0,\infty} = \left\| \widetilde{\xi}_t^u (\mathrm{R}_\tau v_{0,\infty}) - \widetilde{\xi}_t^u (\mathrm{R}_\tau w_{0,\infty}) \right\|_{t,\infty}$$

$$= \left\| F \left( u^t \!\!\leftarrow\!\!\!\rightarrow \mathrm{R}_\tau v_{0,\infty} \right) - F \left( u^t \!\!\leftarrow\!\!\!\rightarrow \mathrm{R}_\tau w_{0,\infty} \right) \right\|_{t,\infty}$$

$$\leq \left\| F \left( u^t \!\!\leftarrow\!\!\!\rightarrow \mathrm{R}_\tau v_{0,\infty} \right) - F \left( u^t \!\!\leftarrow\!\!\!\rightarrow \mathrm{R}_\tau w_{0,\infty} \right) \right\| < \varepsilon \ ,$$

provided

$$\left\| u^t \!\!\leftarrow\!\!\!\rightarrow \mathrm{R}_\tau v_{0,\infty} - u^t \!\!\leftarrow\!\!\!\rightarrow \mathrm{R}_\tau w_{0,\infty} \right\| \leq \|v_{0,\infty} - w_{0,\infty}\|_{t,\infty} < \delta$$

for $\delta = \delta(\varepsilon)$ sufficiently small. Uniform continuity is also inherited by the natural state from the input-output system.

# APPENDIX B

# PROOF OF PROPOSITION 20

## SCALAR CASE

In this appendix, when two systems are considered, one will be denoted by $\xi_t^u(v_{0,\infty})$ and the other will be denoted by $\eta_t^w(v_{0,\infty})$. Differences between $u$ and $w$ will be represented by $z$. At times, $u$ and $w$ will have subscripts. Even when this is the case, subscripts will usually be omitted from differences represented by $z$.

**Integral Operator Case**

We propose to consider the $N-$degree scalar polynomial case of our problem. First, we consider a single $N-$degree term. The system is of the form: (symmetric kernel is assumed)

$$[F(u)](s) = \int_0^\infty \cdots \int_0^\infty f(\tau_1, \cdots, \tau_N) u(s-\tau_1) \cdots u(s-\tau_N) d\tau_1 \cdots d\tau_N \ .$$

The state is of the form:

$$[\xi_0^u(v_{0,\infty})](\sigma) = [F(u_0 \longmapsto v_{0,\infty})](\sigma)$$

$$= \int_0^\infty \cdots \int_0^\infty f(\tau_1, \cdots, \tau_N)(u_0 \longmapsto v_{0,\infty})(\sigma - \tau_1)$$
$$\cdots (u_0 \longmapsto v_{0,\infty})(\sigma - \tau_N) d\tau_1, \cdots, d\tau_N$$

$$= \binom{N}{0} \int_\sigma^\infty \cdots \int_\sigma^\infty f(\tau_1, \cdots, \tau_N) u(\sigma - \tau_1) \cdots u(\sigma - \tau_N) d\tau_1 \cdots d\tau_N$$
$$+ \binom{N}{1} \int_0^\sigma \int_\sigma^\infty \cdots \int_\sigma^\infty f(\tau_1, \cdots, \tau_N) v(\sigma - \tau_1) u(\sigma - \tau_2) \cdots u(\sigma - \tau_N) d\tau_1 \cdots d\tau_N$$
$$+ \binom{N}{2} \int_0^\sigma \int_0^\sigma \int_\sigma^\infty \cdots \int_\sigma^\infty f(\tau_1, \cdots, \tau_N) v(\sigma - \tau_1) v(\sigma - \tau_2) u(\sigma - \tau_3) \cdots u(\sigma - \tau_N) d\tau_1 \cdots d\tau_N$$
$$+ \cdots$$
$$+ \binom{N}{N} \int_0^\sigma \cdots \int_0^\sigma f(\tau_1, \cdots, \tau_N) v(\sigma - \tau_1) \cdots v(\sigma - \tau_N) d\tau_1 \cdots d\tau_N$$

where the integral signs and the $d\tau_i^{\prime s}$ match left to right. Now, replace $v$ by $xv$, where $x \in \Re$,

$$= \binom{N}{0} \widehat{0F}\, x^0 + \binom{N}{1} \widehat{1F}\, x^1 + \binom{N}{2} \widehat{2F}\, x^2 + \cdots + \binom{N}{N} \widehat{NF}\, x^N$$

where

$$\widehat{nF} = \underbrace{\int_0^\sigma \cdots \int_0^\sigma}_{n} \underbrace{\int_\sigma^\infty \cdots \int_\sigma^\infty}_{N-n} f(\tau_1, \cdots, \tau_N) \underbrace{v(\sigma - \tau_{n+1}) \cdots v(\sigma - \tau_N)}_{n}$$

$$\underbrace{u(\sigma-\tau_1)\cdots u(\sigma-\tau_n)}_{N-n} d\tau_1\cdots d\tau_N \ .$$

Two systems are being dealt with as before, and as before, put $h = f - g$; also, since the polynomial is zero, each of the coefficients is zero. Hence,

$$\int_0^\sigma \cdots \int_0^\sigma h(\tau_1, \cdots, \tau_N) v(\sigma - \tau_1) \cdots v(\sigma - \tau_N) d\tau_1 \cdots d\tau_N = 0$$

for all $v$, for all $\sigma > 0$. Suppose $h(a_1, a_2, \cdots a_N) > 0$ and the $a_i'^s$ are all different. Set $A_i = [a_i - \delta/2, a_i + \delta/2]$ where $\delta > 0$ is small enough such that $A_i \cap A_j = \emptyset$ for $i \neq j$. Also, $\delta$ is small enough such that $h(\tau_1, \cdots, \tau_N) > 0$ for all $\tau_i \in A_i$. Hence, with $v(s) = 1$ on $A_i$ ; $= 0$ else,

$$\int_{A_i} \cdots \int_{A_i} h(\tau_1, \cdots, \tau_N) d\tau_1 \cdots d\tau_N = 0 \qquad \text{(B-1)}$$

for all $i$. Next, let

$$v(s) = \begin{cases} 1, & s \in A_i \cup A_j \\ 0, & \text{else} \end{cases}$$

$$\int_{A_i \cup A_j} \cdots \int_{A_i \cup A_j} h(\tau_1, \cdots, \tau_N) d\tau_1 \cdots d\tau_N$$

$$= \binom{N}{0} \int_{A_i} \cdots \int_{A_i} h(\tau_1, \cdots, \tau_N) d\tau_1 \cdots d\tau_N$$

$$+ \binom{N}{1} \int_{A_i} \cdots \int_{A_i} \int_{A_j} h(\tau_1, \cdots, \tau_N) d\tau_1 \cdots d\tau_N$$

$$+ \cdots + \binom{N}{N} \int_{A_j} \cdots \int_{A_j} h(\tau_1, \cdots, \tau_N) d\tau_1 \cdots d\tau_N = 0 \ .$$

This with the above gives

$$\binom{N}{1} \int_{A_i} \cdots \int_{A_i} \int_{A_j} h(\tau_1, \cdots, \tau_N) d\tau_1 \cdots d\tau_N$$

$$+ \cdots + \binom{N}{N-1} \int_{A_i} \int_{A_j} \cdots \int_{A_j} h(\tau_1, \cdots, \tau_N) d\tau_1 \cdots d\tau_N = 0 \ .$$

Next,

$$v(s) = \begin{cases} 1, & s \in A_i \cup A_j \cup A_k \\ 0, & \text{else} \end{cases}$$

$$\int\limits_{A_i \cup A_j \cup A_k} \cdots \int\limits_{A_i \cup A_j \cup A_k} h(\tau_1, \cdots, \tau_N) d\tau_1 \cdots d\tau_N$$

$$= \int\limits_{A_i} + \int\limits_{A_j} + \int\limits_{A_k} + \left\{ \int\limits_{A_i \cup A_j} - \int\limits_{A_i} - \int\limits_{A_j} \right\} + \left\{ \int\limits_{A_i \cup A_k} - \int\limits_{A_i} - \int\limits_{A_k} \right\} + \left\{ \int\limits_{A_j \cup A_k} - \int\limits_{A_j} - \int\limits_{A_k} \right\} + \text{III} .$$

This has three types of integrals: Type 1 is where only a single integral interval appears, i.e., only "$A_i$", "$A_j$" or "$A_k$." These are referred to as integrals with $A'^s$ taken one at a time. By the above, these integrals are zero. Type 2 is where all possible integrals with exactly two intervals appear, as $A_i$ and $A_j$ or $A_i$ and $A_k$. These appear in the parentheses above. By the above, these are also zero. The type 3 integrals are represented by the "III" in the above equation. This type consists of all integrals where $A_i$, $A_j$, and $A_k$ all appear. Since the whole thing is zero and the first two types were zero, we conclude that the type three case of integrals is zero.

We proceed to consider

$$\int\limits_{A_1 \cup A_2 \cup \cdots \cup A_N} \cdots \int h(\tau_1, \cdots, \tau_N) d\tau_1 \cdots d\tau_N \ .$$

This integral is split up into the sum of $N$ different classes or types of integrals. The first type is

$$\int\limits_{A_1} \cdots \int\limits_{A_1} h(\tau_1, \cdots, \tau_N) d\tau_1 \cdots d\tau_N + \cdots + \int\limits_{A_N} \cdots \int\limits_{A_N} h(\tau_1, \cdots, \tau_N) d\tau_1 \cdots d\tau_N \ .$$

It has been found that these are zero. The last type is

$$\int\limits_{A_1} \int\limits_{A_2} \cdots \int\limits_{A_N} h(\tau_1, \cdots, \tau_N) d\tau_1 \cdots d\tau_N \ . \tag{B-2}$$

The types 1 through N are integrals with the appropriate number of $A_i'^s$ appearing. As we have seen, each of the lower integral types is zero; hence, the last integral type is zero. Since the kernels $f$ and $g$ are continuous, $h = 0$, so we have $f(\tau_1, \cdots, \tau_N) = g(\tau_1, \cdots, \tau_N) \ \forall \tau_1, \cdots, \tau_N$.

### $Nth$ Order Scalar Polynomial Case

Following, we consider the $Nth$ order nonhomogeneous case. The state is defined by

$$[\xi_0^u(v_{0,\infty})](\sigma) = f_0 + \int_0^\infty f_1(\tau)[u_0 \mapsto v_{0,\infty}](\sigma - \tau) d\tau$$

$$+ \int_0^\infty \int_0^\infty f_2(\tau_1, \tau_2)[u_0 \mapsto v_{0,\infty}](\sigma - \tau_1)[u_0 \mapsto v_{0,\infty}](\sigma - \tau_2) d\tau_1 d\tau_2$$

$$+ \cdots + \int_0^\infty \cdots \int_0^\infty f_N(\tau_1, \cdots, \tau_N)[u_0 \mapsto v_{0,\infty}](\sigma - \tau_1) \cdots [u_0 \mapsto v_{0,\infty}](\sigma - \tau_N) d\tau_1 \cdots d\tau_N$$

$$= \binom{0}{0} f_0 + \binom{1}{0} \int_\sigma^\infty f_1(\tau) u(\sigma - \tau) d\tau + \binom{1}{1} \int_0^\sigma f_1(\tau) v(\sigma - \tau) d\tau$$

$$+\binom{2}{0}\int_\sigma^\infty \int_\sigma^\infty f_2(\tau_1,\tau_2)u(\sigma-\tau_1)u(\sigma-\tau_2)d\tau_1 d\tau_2$$

$$+\binom{2}{1}\int_0^\sigma \int_\sigma^\infty f_2(\tau_1,\tau_2)v(\sigma-\tau_1)u(\sigma-\tau_2)d\tau_1 d\tau_2$$

$$+\binom{2}{2}\int_0^\sigma \int_0^\sigma f_2(\tau_1,\tau_2)v(\sigma-\tau_1)v(\sigma-\tau_2)d\tau_1 d\tau_2$$

$$+\cdots+\binom{N}{0}\int_\sigma^\infty \cdots \int_\sigma^\infty f_N(\tau_1,\cdots,\tau_N)u(\sigma-\tau_1)\cdots u(\sigma-\tau_N)d\tau_1\cdots d\tau_N$$

$$+\binom{N}{1}\int_0^\sigma \int_\sigma^\infty \cdots \int_\sigma^\infty f_N(\tau_1,\cdots,\tau_N)v(\sigma-\tau_1)u(\sigma-\tau_2)\cdots u(\sigma-\tau_N)d\tau_1\cdots d\tau_N$$

$$+\cdots+\binom{N}{N}\int_0^\sigma \cdots \int_0^\sigma f_N(\tau_1,\cdots,\tau_N)v(\sigma-\tau_1)\cdots v(\sigma-\tau_N)d\tau_1\cdots d\tau_N\ .$$

Now collecting the terms by degree in $v$:

$$=\binom{0}{0}f_0+\binom{1}{0}\int_\sigma^\infty f_1(\tau)u(\sigma-\tau)d\tau+\binom{2}{0}\int_\sigma^\infty \int_\sigma^\infty f_2(\tau_1,\tau_2)u(\sigma-\tau_1)u(\sigma-\tau_2)d\tau_1 d\tau_2$$

$$+\cdots+\binom{N}{0}\int_\sigma^\infty \cdots \int_\sigma^\infty f_N(\tau_1,\cdots,\tau_N)u(\sigma-\tau_1)\cdots u(\sigma-\tau_N)d\tau_1\cdots d\tau_N$$

$$+\binom{1}{1}\int_0^\sigma f_1(\tau)v(\sigma-\tau)d\tau+\binom{2}{1}\int_0^\sigma \int_\sigma^\infty f_2(\tau_1,\tau_2)v(\sigma-\tau_1)u(\sigma-\tau_2)d\tau_1 d\tau_2$$

$$+\cdots+\binom{N}{1}\int_0^\sigma \int_\sigma^\infty \cdots \int_\sigma^\infty f_N(\tau_1,\cdots,\tau_N)v(\sigma-\tau_1)u(\sigma-\tau_2)\cdots u(\sigma-\tau_N)d\tau_1\cdots d\tau_N$$

$$+\cdots+\binom{N-1}{N-1}\int_0^\sigma \cdots \int_0^\sigma f_{N-1}(\tau_1,\cdots,\tau_{N-1})v(\sigma-\tau_1)\cdots v(\sigma-\tau_{N-1})d\tau_1\cdots d\tau_{N-1}$$

$$+\binom{N}{N-1}\int_0^\sigma \cdots \int_0^\sigma \int_\sigma^\infty f_N(\tau_1,\cdots,\tau_N)v(\sigma-\tau_1)\cdots v(\sigma-\tau_{N-1})u(\sigma-\tau_N)d\tau_1\cdots d\tau_N$$

$$+\binom{N}{N}\int_0^\sigma \cdots \int_0^\sigma f_N(\tau_1,\cdots,\tau_N)v(\sigma-\tau_1)\cdots v(\sigma-\tau_N)d\tau_1\cdots d\tau_N\ .$$

Considering two systems with kernels $f$ and $g$, and with the difference of these having kernel $h$, and using the polynomial argument again, each of the terms by degree in $v$ are individually zero. By what has previously been established, for the $N-$degree (in $v$) term we get

$$\int_0^\sigma \cdots \int_0^\sigma h_N(\tau_1,\cdots,\tau_N)v(\sigma-\tau_1)\cdots v(\sigma-\tau_N)d\tau_1\cdots d\tau_N=0$$

and so $h_N(\tau_1, \cdots, \tau_N) = 0$ for all $\tau_1, \cdots \tau_N$. Consider the $(N-1)-$degree (in $v$) terms:

$$\int_0^\sigma \cdots \int_0^\sigma \left\{ h_{N-1}(\tau_1, \cdots, \tau_{N-1}) + \binom{N}{N-1} \int_\sigma^\infty f_N(\tau_1, \cdots, \tau_N) u(\sigma - \tau_N) d\tau_N \right.$$

$$\left. - \binom{N}{N-1} \int_\sigma^\infty g_N(\tau_1, \cdots \tau_N) w\left(\sigma - \tau_N\right) d\tau_N \right\}$$

$$v(\sigma - \tau_1) \cdots v(\sigma - \tau_{N-1}) d\tau_1 \cdots d\tau_{N-1} = 0 \ .$$

If the terms inside the parentheses are continuous, and since $f_N = g_N$, setting $z = w - u$

$$h_{N-1}(\tau_1, \cdots, \tau_{N-1}) = -\binom{N}{N-1} \int_\sigma^\infty f_N(\tau_1, \cdots, \tau_N) z(\sigma - \tau_N) d\tau_N \ .$$

Continuing,

$$|h_{N-1}(\tau_1, \cdots, \tau_{N-1})| \leq C \int_\sigma^\infty |f_N(\tau_1, \cdots, \tau_N)| \, d\tau_N \ ,$$

where $C$ is some positive constant ($z$ is assumed bounded). Taking the limit as $\sigma \to \infty$ yields $h_{N-1}(\tau_1, \cdots, \tau_{N-1}) = 0$ for all $\tau_1, \cdots, \tau_{N-1}$, since we assume that $f_N$ is absolutely integrable. Now going back through the terms in decreasing order,

$$\int_0^\sigma \left[ h_1(\tau_1) + \left\{ \binom{2}{1} \int_\sigma^\infty f_2(\tau_1, \tau_2) z_2(\sigma - \tau_2) d\tau_2 + \cdots \right. \right.$$

$$\left. \left. + \binom{N}{1} \int_\sigma^\infty \cdots \int_\sigma^\infty f_N(\tau_1, \cdots, \tau_N) z_N(\sigma - \tau_2, \cdots, \sigma - \tau_N) d\tau_2 \cdots d\tau_N \right\} \right] v(\sigma - \tau_1) d\tau_1 = 0.$$

If the terms inside the parentheses are continuous, and $u$, $w$ are bounded, then

$$|h_1(\tau_1)| \leq C_2 \int_\sigma^\infty |f_2(\tau_1, \tau_2)| \, d\tau_2 + \cdots + C_N \int_\sigma^\infty \cdots \int_\sigma^\infty |f_N(\tau_1, \cdots, \tau_N)| \, d\tau_2 \cdots d\tau_N \ ,$$

where the $C_i'^s$ are positive constants. Taking the limit as $\sigma \to \infty$ gives $h_1(\tau) = 0$ for all $\tau$. Looking at the zero<u>th</u> degree term,

$$|h_0| \leq C_1 \int_\sigma^\infty |f_1(\tau)| \, d\tau + \cdots + C_N \int_\sigma^\infty \cdots \int_\sigma^\infty |f_N(\tau_1, \cdots, \tau_N)| \, d\tau_1 \cdots d\tau_N \ ,$$

where the $C_i'^s$ are positive constants. (The $C_i'^s$ in separate equations are not the same.) Taking the limit $\sigma \to \infty$ gives $h_0 = 0$.

We summarize the necessary conditions. All the kernels $f_1, \cdots f_N$ are absolutely integrable and continuous. These terms also need be continuous:

$$\int_\sigma^\infty f_N(\tau_1, \cdots, \tau_N) z(\sigma - \tau_N) d\tau_N \ ;$$

$$\int_\sigma^\infty \int_\sigma^\infty f_N(\tau_1, \cdots \tau_N) z(\sigma - \tau_{N-1}, \sigma - \tau_N) d\tau_{N-1} d\tau_N \ ;$$

$$\vdots$$

$$\int_{\sigma}^{\infty}\cdots\int_{\sigma}^{\infty} f_N(\tau_1,\cdots\tau_N)z(\sigma-\tau_2,\cdots,\sigma-\tau_N)d\tau_2\cdots d\tau_N \ ;$$

$$\int_{\sigma}^{\infty} f_{N-1}(\tau_1,\cdots\tau_{N-1})z(\sigma-\tau_{N-1})d\tau_{N-1} \ ;$$

$$\vdots$$

$$\int_{\sigma}^{\infty}\cdots\int_{\sigma}^{\infty} f_{N-1}(\tau_1,\cdots\tau_{N-1})z(\sigma-\tau_2,\cdots,\sigma-\tau_{N-1})d\tau_2\cdots d\tau_{N-1} \ ;$$

$$\vdots$$

$$\int_{\sigma}^{\infty} f_2(\tau_1,\tau_2)z(\sigma-\tau_2)d\tau_2 \quad .$$

We consider sufficient conditions for the continuity of these terms.

$$\lim_{\substack{\tau_1',\cdots\tau_i' \\ \to\tau_1,\cdots,\tau_i}} \int_{\sigma}^{\infty}\cdots\int_{\sigma}^{\infty} \big|f_N(\tau_1',\cdots,\tau_i',\tau_{i+1},\cdots,\tau_N)z(\sigma-\tau_{i+1},\cdots,\sigma-\tau_N)$$

$$-f_N(\tau_1,\cdots,\tau_i,\tau_{i+1},\cdots,\tau_N)z(\sigma-\tau_{i+1},\cdots,\sigma-\tau_N)\big|\, d\tau_{i+1}\cdots d\tau_N$$

$$\le Const. \lim_{\substack{\tau_1',\cdots,\tau_i' \\ \to\tau_1,\cdots,\tau_i}} \int_{\sigma}^{\infty}\cdots\int_{\sigma}^{\infty} \big|f_N(\tau_1',\cdots,\tau_i',\tau_{i+1},\cdots,\tau_N)$$

$$-f_N(\tau_1,\cdots,\tau_i,\tau_{i+1},\cdots,\tau_N)\big|\, d\tau_{i+1}\cdots d\tau_N \ .$$

Assume

$$\sup_{\tau_1,\cdots,\tau_i} f_N(\tau_1,\cdots,\tau_i,\tau_{i+1},\cdots,\tau_N) \triangleq g(\tau_{i+1},\cdots,\tau_N)$$

is absolutely integrable (as a function of $\tau_{i+1},\cdots,\tau_N$). Hence, from the Lebesgue Dominated Convergence Theorem:

$$\lim_{\substack{\tau_1',\cdots,\tau_i' \\ \to\tau_1,\cdots,\tau_i}} \int_{\sigma}^{\infty}\cdots\int_{\sigma}^{\infty} \big|f_N(\tau_1',\cdots,\tau_i',\tau_{i+1},\cdots,\tau_N)-f_N(\tau_1,\cdots,\tau_i,\tau_{i+1},\cdots,\tau_N)\big|\, d\tau_{i+1}\cdots d\tau_N$$

$$=\int_{\sigma}^{\infty}\cdots\int_{\sigma}^{\infty} \lim_{\substack{\tau_1',\cdots,\tau_i' \\ \to\tau_1,\cdots,\tau_i}} \big|f_N(\tau_1',\cdots,\tau_i',\tau_{i+1},\cdots,\tau_N)$$

$$-f_N(\tau_1,\cdots,\tau_i,\tau_{i+1},\cdots,\tau_N)\big|\, d\tau_{i+1}\cdots d\tau_N = 0 \ .$$

## VECTOR CASE

### Definition of Operator and Permutation

The system at hand is an $N-$degree homogeneous integral operator with $M$-dimensional input space. An input $u$ is of the form $[u_1, \cdots, u_M]$ with norm $\|u\|_{s,t} = \max_i \|u_i\|_{s,t}$ and $\|u_i\|_{s,t} = \sup_{s<\tau\leq t} |u_i(\tau)|$. Without loss of generality, the output space has dimension one. The system mapping is represented by

$$[F(u)](t) = \int_0^\infty \cdots \int_0^\infty \sum_{i_1,\cdots,i_N=1}^{M} f_N^{i_1,\cdots,i_N}(\sigma_1, \cdots, \sigma_N) u_{i_1}(t-\sigma_1) \cdots u_{i_N}(t-\sigma_N) d\sigma_1 \cdots d\sigma_N \ ,$$

where $i_1, \cdots, i_N \in \{1, 2, \cdots, M\}$, not necessarily distinct. Note that $M$ may be greater than, less than, or equal to $N$. For example, if $N = 2$ and $M = 3$, $[i_1, i_2] = (1,1)$, $(2,2)$, $(3,3)$, $(1,2)$, $(1,3)$, $(2,1)$, $(2,3)$, $(3,1)$, $(3,2)$. That is, these are the values $[i_1, i_2]$ take on in the above sum. If $N = 3$ and $M = 2$, $[i_1, i_2, i_3] = (1,1,1)$, $(2,2,2)$, $(1,1,2)$, $(1,2,1)$, $(2,1,1)$, $(1,2,2)$, $(2,1,2)$, $(2,2,1)$. There are $M^N$ terms in this sum.

A proposed symmetrization for each kernel $f$ is

$$\widetilde{f}_N^{i_1,\cdots,i_N}(t_1, \cdots, t_N) = \frac{1}{N!} \sum_{\pi} f_N^{\pi(i_1),\cdots,\pi(i_N)}(t_{\pi(1)}, \cdots, t_{\pi(N)})$$

where $\pi$ is the permutation of $N$ ordered symbols; if $\pi(a_1, \cdots, a_N) = (b_1, \cdots, b_N)$ we write $\pi(a_i) = b_i$. For example, in the above sum, if $N = 3$ and $[i_1, i_2, i_3] = (1,2,3)$, the terms appearing would have $(1,2,3)$, $(1,3,2)$, $(2,1,3)$, $(2,3,1)$, $(3,1,2)$, $(3,2,1)$ for the superscripts on $f_N$. If $N = 3$ and $[i_1, i_2, i_3] = (1,3,3)$, the terms appearing would be $(1,3,3)$, $(1,3,3)$, $(3,1,3)$, $(3,1,3)$, $(3,3,1)$, $(3,3,1)$. There are $N!$ terms in this sum. Note that we also have that

$$\widetilde{f}_N^{i_1,\cdots,i_N}(t_1, \cdots, t_N) = \frac{1}{N!} \sum_{\pi} f_N^{i_{\pi(1)},\cdots,i_{\pi(N)}}(t_{\pi(1)}, \cdots, t_{\pi(N)})$$

and we will use this form. This symmetrization has the following two properties:

(i) - For the first property consider,

$$\int_0^\infty \cdots \int_0^\infty \sum_{i_1,\cdots,i_N=1}^{M} \widetilde{f}_N^{i_1,\cdots,i_N}(\sigma_1, \cdots, \sigma_N) u_{i_1}(t-\sigma_1) \cdots u_{i_N}(t-\sigma_N) d\sigma_1 \cdots d\sigma_N$$

$$= \int_0^\infty \cdots \int_0^\infty \sum_{i_1,\cdots,i_N=1}^{M} \left\{ \frac{1}{N!} \sum_{\pi} f_N^{i_{\pi(1)},\cdots,i_{\pi(N)}}(\sigma_{\pi(1)}, \cdots, \sigma_{\pi(N)}) \right\}$$

$$u_{i_1}(t-\sigma_1) \cdots u_{i_N}(t-\sigma_N) d\sigma_1 \cdots d\sigma_N \ . \qquad \text{(B-3)}$$

Now, consider a specific set of indices $i_1, \cdots, i_N$. Suppose that there are $I$ distinct indices and a particular index $i_j$ is repeated $r_j$ times. Then $r_1 + r_2 + \cdots + r_I = N$. We will say that the outside sum cares about repeats (takes repeats into account), but the inside sum does not take repeats into account. Between the two sums, a total of $N! \cdot \binom{N}{r_1}\binom{N-r_1}{r_2} \cdots \binom{N-(r_1+\cdots+r_{I-1})}{r_I}$ number of terms are generated for the specific set under consideration. But, there are only $\binom{N}{r_1}\binom{N-r_1}{r_2}\cdots\binom{N-(r_1+\cdots+r_{I-1})}{r_I}$ different ways to permute this specific set. In each term, there are four sets of indices:

1. The "superscript" indices on "$f$"
2. The indices associated with the $\sigma'^s$ on the argument of $f$
3. The indices on the $u'^s$
4. The indices on the arguments of the $u'^s$

First, assume no repeats in the outer sum, i.e., all $r_i = 1$. Note that $\binom{N}{1}\binom{N-1}{1}\cdots\binom{2}{1}\binom{1}{1} = N!$. With the inner sum, in this case, there are $N!\cdot N!$ terms. Each of the $N!$ outside sum "base" terms (a term on the left-hand side of (B-3)), is permuted $N!$ times by the inner sum. Take a specific outside "base" term with the indices in a specific order. Of the $N!$ terms it gives, only the identity permutation term has its indices in the specified order:

$$f^{1,\cdots,N}(\sigma_1,\cdots,\sigma_N)u_1(t-\sigma_1)\cdots u_N(t-\sigma_N)\ . \tag{B-4}$$

However, each of the other outside "base" terms, of the form

$$\tilde{f}^{\pi(1),\cdots,\pi(N)}(\sigma_1,\cdots,\sigma_N)u_{\pi(1)}(t-\sigma_1)\cdots u_{\pi(N)}(t-\sigma_N)\ ,$$

has a single term with the permutation $\pi^{-1}$, which has indices in the specified order

$$f^{1,\cdots,N}(\sigma_{\pi^{-1}(1)},\cdots,\sigma_{\pi^{-1}(N)})u_{\pi(1)}(t-\sigma_1)\cdots u_{\pi(N)}(t-\sigma_N)$$

$$= f^{1,\cdots,N}(\sigma_{\pi^{-1}(1)},\cdots,\sigma_{\pi^{-1}(N)})u_1(t-\sigma_{\pi^{-1}(1)})\cdots u_N(t-\sigma_{\pi^{-1}(N)})\ . \tag{B-5}$$

However, (B-4) and (B-5) are identical, since there are $N!$ such terms and since the inside part of (B-3) is divided by $N!$, we see that under the condition of all $r_i = 1$ the symmetrization does not alter the operator.

Next, suppose $r_1 \neq 1$ and all other $r_i = 1$. With the indices in their specific order, the outside term is

$$\tilde{f}^{1,\cdots 1,2,\cdots,P}(\sigma_1,\sigma_2,\cdots,\sigma_N)u_1(t-\sigma_1)\cdots u_1(t-\sigma_{r_1})u_2(t-\sigma_{r_1+1})\cdots u_P(t-\sigma_N)\ .$$

This term gives $r_1!$ inner terms with indices in this order, and they are all identical. Now, there are $\binom{N}{r_1}\cdot(N-r_1)!$ "outside" terms with these indices in any order; so far we have considered one, the one with the indices in the specified order. Each of these other "outside" terms gives $r_1!$ inner terms with these indices in the specified order. Total number of terms in the specified order are $\binom{N}{r_1}\cdot(N-r_1)!\cdot r_1! = N!$. However, dividing by $N!$ gives the symmetrization and once again does not alter the operator.

Now for a couple of examples: First, consider three distinct indices

$$\tilde{f}^{123}(\sigma_1,\sigma_2,\sigma_3) = \frac{1}{3!}[f^{123}(\sigma_1,\sigma_2,\sigma_3) + f^{213}(\sigma_2,\sigma_1,\sigma_3) + f^{312}(\sigma_3,\sigma_1,\sigma_2)$$

$$+f^{132}(\sigma_1,\sigma_3,\sigma_2) + f^{231}(\sigma_2,\sigma_3,\sigma_1) + f^{321}(\sigma_3,\sigma_2,\sigma_1)$$

$$\tilde{f}^{213}(\sigma_1,\sigma_2,\sigma_3) = \frac{1}{3!}[f^{123}(\sigma_2,\sigma_1,\sigma_3) + f^{213}(\sigma_1,\sigma_2,\sigma_3) + f^{312}(\sigma_3,\sigma_2,\sigma_1)$$
$$+f^{132}(\sigma_2,\sigma_3,\sigma_1) + f^{231}(\sigma_1,\sigma_3,\sigma_2) + f^{321}(\sigma_3,\sigma_1,\sigma_2)$$
$$\tilde{f}^{312}(\sigma_1,\sigma_2,\sigma_3) = \frac{1}{3!}[f^{123}(\sigma_2,\sigma_3,\sigma_1) + f^{213}(\sigma_3,\sigma_2,\sigma_1) + f^{312}(\sigma_1,\sigma_2,\sigma_3)$$
$$+f^{132}(\sigma_2,\sigma_1,\sigma_3) + f^{231}(\sigma_3,\sigma_1,\sigma_2) + f^{321}(\sigma_1,\sigma_3,\sigma_2)$$
$$\tilde{f}^{132}(\sigma_1,\sigma_2,\sigma_3) = \frac{1}{3!}[f^{123}(\sigma_1,\sigma_3,\sigma_2) + f^{213}(\sigma_3,\sigma_1,\sigma_2) + f^{312}(\sigma_2,\sigma_1,\sigma_3)$$
$$+f^{132}(\sigma_1,\sigma_2,\sigma_3) + f^{231}(\sigma_3,\sigma_2,\sigma_1) + f^{321}(\sigma_2,\sigma_3,\sigma_1)$$
$$\tilde{f}^{231}(\sigma_1,\sigma_2,\sigma_3) = \frac{1}{3!}[f^{123}(\sigma_3,\sigma_1,\sigma_2) + f^{213}(\sigma_1,\sigma_3,\sigma_2) + f^{312}(\sigma_2,\sigma_3,\sigma_1)$$
$$+f^{132}(\sigma_3,\sigma_2,\sigma_1) + f^{231}(\sigma_1,\sigma_2,\sigma_3) + f^{321}(\sigma_2,\sigma_1,\sigma_3)$$
$$\tilde{f}^{321}(\sigma_1,\sigma_2,\sigma_3) = \frac{1}{3!}[f^{123}(\sigma_3,\sigma_2,\sigma_1) + f^{213}(\sigma_2,\sigma_3,\sigma_1) + f^{312}(\sigma_1,\sigma_3,\sigma_2)$$
$$+f^{132}(\sigma_3,\sigma_1,\sigma_2) + f^{231}(\sigma_2,\sigma_1,\sigma_3) + f^{321}(\sigma_1,\sigma_2,\sigma_3)\ .$$

Sum and integrate the $f^{231}$ terms. For example:

$$\int f^{231}(\sigma_2,\sigma_3,\sigma_1)u_1(t-\sigma_1)u_2(t-\sigma_2)u_3(t-\sigma_3)$$
$$+\int f^{231}(\sigma_1,\sigma_3,\sigma_2)u_2(t-\sigma_1)u_1(t-\sigma_2)u_3(t-\sigma_3)$$
$$+\int f^{231}(\sigma_3,\sigma_1,\sigma_2)u_3(t-\sigma_1)u_1(t-\sigma_2)u_2(t-\sigma_3)$$
$$+\int f^{231}(\sigma_3,\sigma_2,\sigma_1)u_1(t-\sigma_1)u_3(t-\sigma_2)u_2(t-\sigma_3)$$
$$+\int f^{231}(\sigma_1,\sigma_2,\sigma_3)u_2(t-\sigma_1)u_3(t-\sigma_2)u_1(t-\sigma_3)$$
$$+\int f^{231}(\sigma_2,\sigma_1,\sigma_3)u_3(t-\sigma_1)u_2(t-\sigma_2)u_1(t-\sigma_3)\ .$$

These terms are identical. For the second example, consider the $N = 3$, $M = 2$ case given above. In this example, the left-hand side of (B-3) is

$$\int \tilde{f}^{111}(\sigma_1,\sigma_2,\sigma_3)u_1(t-\sigma_1)u_1(t-\sigma_2)u_1(t-\sigma_3)$$
$$+\int \tilde{f}^{112}(\sigma_1,\sigma_2,\sigma_3)u_1(t-\sigma_1)u_1(t-\sigma_2)u_2(t-\sigma_3)$$

$$+\int \tilde{f}^{121}(\sigma_1,\sigma_2,\sigma_3)u_1(t-\sigma_1)u_2(t-\sigma_2)u_1(t-\sigma_3)$$

$$+\int \tilde{f}^{122}(\sigma_1,\sigma_2,\sigma_3)u_1(t-\sigma_1)u_2(t-\sigma_2)u_2(t-\sigma_3)$$

$$+\int \tilde{f}^{211}(\sigma_1,\sigma_2,\sigma_3)u_2(t-\sigma_1)u_1(t-\sigma_2)u_1(t-\sigma_3)$$

$$+\int \tilde{f}^{212}(\sigma_1,\sigma_2,\sigma_3)u_2(t-\sigma_1)u_1(t-\sigma_2)u_2(t-\sigma_3)$$

$$+\int \tilde{f}^{221}(\sigma_1,\sigma_2,\sigma_3)u_2(t-\sigma_1)u_2(t-\sigma_2)u_1(t-\sigma_3)$$

$$+\int \tilde{f}^{222}(\sigma_1,\sigma_2,\sigma_3)u_2(t-\sigma_1)u_2(t-\sigma_2)u_2(t-\sigma_3)\ .$$

Take all terms with $r_1 = 2$ and $r_2 = 1$. Expand these using the inner sum,

$$\tilde{f}^{112}(\sigma_1,\sigma_2,\sigma_3) = \frac{1}{3!}\{f^{112}(\sigma_1,\sigma_2,\sigma_3) + f^{112}(\sigma_2,\sigma_1,\sigma_3) + f^{121}(\sigma_1,\sigma_3,\sigma_2)$$
$$f^{121}(\sigma_2,\sigma_3,\sigma_1) + f^{211}(\sigma_3,\sigma_1,\sigma_2) + f^{211}(\sigma_3,\sigma_2,\sigma_1)\}$$
$$\tilde{f}^{121}(\sigma_1,\sigma_2,\sigma_3) = \frac{1}{3!}\{f^{121}(\sigma_1,\sigma_2,\sigma_3) + f^{121}(\sigma_3,\sigma_2,\sigma_1) + f^{112}(\sigma_1,\sigma_3,\sigma_2)$$
$$f^{112}(\sigma_3,\sigma_1,\sigma_2) + f^{211}(\sigma_2,\sigma_1,\sigma_3) + f^{211}(\sigma_2,\sigma_3,\sigma_1)\}$$
$$\tilde{f}^{211}(\sigma_1,\sigma_2,\sigma_3) = \frac{1}{3!}\{f^{211}(\sigma_1,\sigma_2,\sigma_3) + f^{211}(\sigma_1,\sigma_3,\sigma_2) + f^{121}(\sigma_2,\sigma_1,\sigma_3)$$
$$f^{121}(\sigma_3,\sigma_1,\sigma_2) + f^{112}(\sigma_2,\sigma_3,\sigma_1) + f^{112}(\sigma_3,\sigma_2,\sigma_1)\}\ .$$

Pull out and put into integrals the terms with indices in the order "1, 1, 2," for example,

$$\frac{1}{3!}\int\int\int[\{f^{112}(\sigma_1,\sigma_2,\sigma_3) + f^{112}(\sigma_2,\sigma_1,\sigma_3)\}u_1(t-\sigma_1)u_1(t-\sigma_2)u_2(t-\sigma_3)$$
$$\{f^{112}(\sigma_1,\sigma_3,\sigma_2) + f^{112}(\sigma_3,\sigma_1,\sigma_2)\}u_1(t-\sigma_1)u_2(t-\sigma_2)u_1(t-\sigma_3)$$
$$\{f^{112}(\sigma_2,\sigma_3,\sigma_1) + f^{112}(\sigma_3,\sigma_2,\sigma_1)\}u_2(t-\sigma_1)u_1(t-\sigma_2)u_1(t-\sigma_3)]d\sigma_1 d\sigma_2 d\sigma_3\ .$$

These terms are all the same. The first property is that symmetrization does not alter the operator.

(ii) For the second property consider,

$$\widetilde{f}_N^{i_1,\cdots,i_N}(\sigma_1,\cdots,\sigma_N)u_{i_1}(t-\sigma_1)\cdots u_{i_N}(t-\sigma_N)$$

$$= \frac{1}{N!}\sum_{\pi} f_N^{i_{\pi(1)},\cdots,i_{\pi(N)}}(\sigma_{\pi(1)},\cdots,\sigma_{\pi(N)})u_{i_1}(t-\sigma_1)\cdots u_{i_N}(t-\sigma_N)\ . \qquad \text{(B-6)}$$

Also

$$\widetilde{f}_N^{i_{\rho(1)},\cdots,i_{\rho(N)}}(\sigma_1,\cdots,\sigma_N)u_{i_{\rho(1)}}(t-\sigma_1)\cdots u_{i_{\rho(N)}}(t-\sigma_N)$$

$$=\frac{1}{N!}\sum_{\pi}f_N^{i_{\pi(\rho(1))},\cdots,i_{\pi(\rho(N))}}(\sigma_{\pi(1)},\cdots,\sigma_{\pi(N)})u_{i_{\rho(1)}}(t-\sigma_1)\cdots u_{i_{\rho(N)}}(t-\sigma_N)\ . \tag{B-7}$$

Show (B-6) equals (B-7) (this is the second property). Take a particular term out of (B-6), that is, fix $\pi$,

$$f_N^{i_{\pi(1)},\cdots,i_{\pi(N)}}(\sigma_{\pi(1)},\cdots,\sigma_{\pi(N)})u_{i_1}(t-\sigma_1)\cdots u_{i_N}(t-\sigma_N)\ .$$

There exists $\pi'$ such that $\pi=\pi'\rho$. Take this term out of (B-7)

$$f_N^{i_{\pi'\rho(1)},\cdots,i_{\pi'\rho(N)}}(\sigma_{\pi'(1)},\cdots,\sigma_{\pi'(N)})u_{i_{\rho(1)}}(t-\sigma_1)\cdots u_{i_{\rho(N)}}(t-\sigma_N)$$

$$=f_N^{i_{\pi(1)},\cdots,i_{\pi(N)}}(\sigma_{\pi'\rho\rho^{-1}(1)},\cdots,\sigma_{\pi'\rho\rho^{-1}(N)})u_{i_{\rho(1)}}(t-\sigma_1)\cdots u_{i_{\rho(N)}}(t-\sigma_N)$$

$$=f_N^{i_{\pi(1)},\cdots,i_{\pi(N)}}(\sigma_{\pi\rho^{-1}(1)},\cdots,\sigma_{\pi\rho^{-1}(N)})u_{i_{\rho(1)}}(t-\sigma_1)\cdots u_{i_{\rho(N)}}(t-\sigma_N)\ .$$

Let $\tau_{\rho(1)}=\sigma_1,\cdots,\tau_{\rho(N)}=\sigma_N$, then $\sigma_{\pi\rho^{-1}(1)}=\tau_{\pi(1)},\cdots,\sigma_{\pi\rho^{-1}(N)}=\tau_{\pi(N)}$. Hence,

$$=f_N^{i_{\pi(1)},\cdots,i_{\pi(N)}}(\tau_{\pi(1)},\cdots,\tau_{\pi(N)})u_{i_{\rho(1)}}(t-\tau_{\rho(1)})\cdots u_{i_{\rho(N)}}(t-\tau_{\rho(N)})\ .$$

Since $\rho$ is 1 to 1,

$$=f_N^{i_{\pi(1)},\cdots,i_{\pi(N)}}(\tau_{\pi(1)},\cdots,\tau_{\pi(N)})u_{i_1}(t-\tau_1)\cdots u_{i_N}(t-\tau_N)\ .$$

As an example, let $N=6$ and $\{i_1,i_2,i_3,i_4,i_5,i_6\}=\{1,2,3,4,5,6\}$. Let the permutations be

$$\begin{array}{rcccccc} & 1 & 2 & 3 & 4 & 5 & 6\\ \pi & 2 & 1 & 3 & 5 & 4 & 6\\ \rho & 3 & 4 & 6 & 5 & 2 & 1\\ \pi^{-1} & 2 & 1 & 3 & 5 & 4 & 6\\ \rho^{-1} & 6 & 5 & 1 & 2 & 4 & 3\\ \pi\rho^{-1} & 5 & 6 & 1 & 4 & 2 & 3 \end{array}\ .$$

Then, using $\tau_{\rho(i)}=\sigma_i$, $\sigma_1=\tau_3$, $\sigma_2=\tau_4$, $\sigma_3=\tau_6$, $\sigma_4=\tau_5$, $\sigma_5=\tau_2$, $\sigma_6=\tau_1$. Computing,

$$\begin{array}{lllll} \tau_{\pi(1)}=\tau_2 & ; & \tau_{\pi(3)}=\tau_3 & ; & \tau_{\pi(5)}=\tau_4\\ \tau_{\pi(2)}=\tau_1 & ; & \tau_{\pi(4)}=\tau_5 & ; & \tau_{\pi(6)}=\tau_6 \end{array}\ .$$

Then,

$$\begin{array}{ll} \sigma_{\pi\rho^{-1}(1)}=\sigma_5=\tau_2=\tau_{\pi(1)} & \sigma_{\pi\rho^{-1}(4)}=\sigma_4=\tau_5=\tau_{\pi(4)}\\ \sigma_{\pi\rho^{-1}(2)}=\sigma_6=\tau_1=\tau_{\pi(2)} & \sigma_{\pi\rho^{-1}(5)}=\sigma_2=\tau_4=\tau_{\pi(5)}\\ \sigma_{\pi\rho^{-1}(3)}=\sigma_1=\tau_3=\tau_{\pi(3)} & \sigma_{\pi\rho^{-1}(6)}=\sigma_3=\tau_6=\tau_{\pi(6)} \end{array}\ .$$

As a second example, try a different $\pi$ with the same $\rho$:

$$\begin{array}{ccccccc} & 1 & 2 & 3 & 4 & 5 & 6 \\ \pi & 4 & 2 & 1 & 5 & 3 & 6 \\ \rho & 3 & 4 & 6 & 5 & 2 & 1 \\ \pi^{-1} & 3 & 2 & 5 & 1 & 4 & 6 \\ \rho^{-1} & 6 & 5 & 1 & 2 & 4 & 3 \\ \pi\rho^{-1} & 2 & 5 & 6 & 4 & 1 & 3 \end{array} .$$

Computing,

$$\begin{array}{ccccc} \tau_{\pi(1)} = \tau_4 & ; & \tau_{\pi(3)} = \tau_1 & ; & \tau_{\pi(5)} = \tau_3 \\ \tau_{\pi(2)} = \tau_2 & ; & \tau_{\pi(4)} = \tau_5 & ; & \tau_{\pi(6)} = \tau_6 \end{array} .$$

Then,

$$\begin{array}{ll} \sigma_{\pi\rho^{-1}(1)} = \sigma_2 = \tau_4 = \tau_{\pi(1)} & \sigma_{\pi\rho^{-1}(4)} = \sigma_4 = \tau_5 = \tau_{\pi(4)} \\ \sigma_{\pi\rho^{-1}(2)} = \sigma_5 = \tau_2 = \tau_{\pi(2)} & \sigma_{\pi\rho^{-1}(5)} = \sigma_1 = \tau_3 = \tau_{\pi(5)} \\ \sigma_{\pi\rho^{-1}(3)} = \sigma_6 = \tau_1 = \tau_{\pi(3)} & \sigma_{\pi\rho^{-1}(6)} = \sigma_3 = \tau_6 = \tau_{\pi(6)} \end{array} .$$

The first property shows that the proposed formula gives the symmetrization and the second property is the one needed in the development.

**Integral Operator Case**

$$[\xi_0^u(v_{0,\infty})](s) = [F(u_0 \longmapsto v_{0,\infty})](s)$$

$$= \int_0^\infty \cdots \int_0^\infty \sum_{i_1,\cdots,i_N=1}^{M} \widetilde{f}_N^{i_1,\cdots,i_N}(\sigma_1,\cdots,\sigma_N)(u_{i_1} \longmapsto v_{i_1})(s-\sigma_1)$$

$$\cdots (u_{i_N} \longmapsto v_{i_N})(s-\sigma_N) d\sigma_1 \cdots d\sigma_N \ .$$

The second system is

$$[\eta_0^w(v_{0,\infty})](s) = [G(w_0 \longmapsto v_{0,\infty})](s)$$

$$= \int_0^\infty \cdots \int_0^\infty \sum_{i_1,\cdots,i_N=1}^{M} \widetilde{g}_N^{i_1,\cdots,i_N}(\sigma_1,\cdots,\sigma_N)(w_{i_1} \longmapsto v_{i_1})(s-\sigma_1)$$

$$\cdots (w_{i_N} \longmapsto v_{i_N})(s-\sigma_N) d\sigma_1 \cdots d\sigma_N \ .$$

Now,

$$[\xi_0^u(v_{0,\infty})](s) = \int_0^\infty \cdots \int_0^\infty \widetilde{f}_N^{i_1,\cdots,i_1}(\sigma_1,\cdots,\sigma_N)(u_{i_1} \longmapsto v_{i_1})(s-\sigma_1)$$

$$\cdots (u_{i_1} \longmapsto v_{i_1})(s-\sigma_N) d\sigma_1 \cdots d\sigma_N$$

$$+ \cdots + \int_0^\infty \cdots \int_0^\infty \widetilde{f}_N^{i_M,\cdots,i_M}(\sigma_1,\cdots,\sigma_N)(u_{i_M} \longmapsto v_{i_M})(s-\sigma_1)$$

$$\cdots (u_{i_M} \longmapsto v_{i_M})(s-\sigma_N) d\sigma_1 \cdots d\sigma_N \ .$$

Consider a summand with the difference kernel "$\widetilde{h}$"

$$\int_0^s \cdots \int_0^s \widetilde{h}_N^{i_1,\cdots,i_N}(\sigma_1,\cdots,\sigma_N)v_{i_1}(s-\sigma_1)\cdots v_{i_N}(s-\sigma_N)d\sigma_1\cdots d\sigma_N \ .$$

The possible input indices $i_1,\cdots,i_N \in \{1,2,\cdots,M\}$ may be listed in order in $\widetilde{h}_N^{i_1,\cdots,i_N}$ since $\widetilde{h}_N$ is a symmetric kernel. These indices may be repeated, such as

$$i_1 = \cdots = i_{r_1};\ i_{r_1+1} = \cdots = i_{r_1+r_2};$$

$$\cdots;\ i_{r_1+r_2+\cdots+r_{I-1}+1} = \cdots = i_{r_1+r_2+\cdots+r_I}$$

such that $r_1+r_2+\cdots+r_I = N$. Using the "$v \to xv$" argument, we may focus on a single summand: (a summand of highest order; i.e., having the greatest number of $v'^s$)

$$\overbrace{\int_{A_1}\cdots\int}^{r_1}\overbrace{\int_{A_2}\cdots\int}^{r_2}\cdots\overbrace{\int_{A_I}\cdots\int}^{r_I}\widetilde{h}_N^{i_1,\cdots,i_N}(\sigma_1,\cdots,\sigma_N)d\sigma_1\cdots d\sigma_N = 0 \tag{B-8}$$

where the $A_i'^s$ are intervals as in the scalar case. For the partitions $A_1, A_2, \cdots, A_{I-1}$ fixed, we use the technique used in the homogeneous scalar case on the $Ith$ (last) partition to show that

$$\overbrace{\int_{A_1}\cdots\int}^{r_1}\overbrace{\int_{A_2}\cdots\int}^{r_2}\cdots\overbrace{\int_{A_{I,1}\cdots A_{I,r_I}}\cdots\int}^{r_I}\widetilde{h}_N^{i_1,\cdots,i_N}(\sigma_1,\cdots,\sigma_N)d\sigma_1\cdots d\sigma_N = 0\ . \tag{B-9}$$

We move to the $I-1th$ partition. We repeat the same procedure. Finally, we find

$$\overbrace{\int_{A_{1,1}\cdots A_{1,r_1}}\cdots\int}^{r_1}\overbrace{\int_{A_{2,1}\cdots A_{2,r_2}}\cdots\int}^{r_2}\cdots\overbrace{\int_{A_{I,1}\cdots A_{I,rI}}\cdots\int}^{r_I}\widetilde{h}_N^{i_1,\cdots,i_N}(\sigma_1,\cdots,\sigma_N)d\sigma_1\cdots d\sigma_N = 0\ .$$

Going from (B-8) to (B-9) is analogous to going from (B-1) to (B-2) in the scalar case. Since $\widetilde{f}$ and $\widetilde{g}$ are continuous, $\widetilde{h}=0$, so we have $\widetilde{f}=\widetilde{g}$.

### *Nth* Order Vector Polynomial Case

$$[F(u_0 \longmapsto v_{0,\infty})](s) = \int_0^\infty \cdots \int_0^\infty \sum_{i_1,\cdots,i_N=1}^{M} \widetilde{f}_N^{i_1,\cdots,i_N}(\sigma_1,\cdots,\sigma_N)(u_{i_1} \longmapsto v_{i_1})(s-\sigma_1)$$

$$\cdots (u_{i_N} \longmapsto v_{i_N})(s-\sigma_N)d\sigma_1\cdots d\sigma_N$$

$$+\int_0^\infty \cdots \int_0^\infty \sum_{i_1,\cdots,i_{N-1}=1}^{M} \widetilde{f}_N^{i_1,\cdots,i_{N-1}}(\sigma_1,\cdots,\sigma_{N-1})(u_{i_1} \longmapsto v_{i_1})(s-\sigma_1)$$

$$\cdots (u_{i_{N-1}} \longmapsto v_{i_{N-1}})(s-\sigma_{N-1})d\sigma_1 \cdots d\sigma_{N-1}$$

$$+\cdots$$

$$+\int_0^\infty \int_0^\infty \sum_{i_1,i_2=1}^{M} \widetilde{f}_N^{i_1,i_2}(\sigma_1,\sigma_2)(u_{i_1} \longmapsto v_{i_1})(s-\sigma_1)(u_{i_2} \longmapsto v_{i_2})(s-\sigma_2)d\sigma_1 d\sigma_2$$

$$+\int_0^\infty \sum_{i=1}^{M} \widetilde{f}_N^{i}(\sigma)(u_i \longmapsto v_i)(s-\sigma)d\sigma + f_0 \ .$$

As with the scalar case, we start with the higher order terms and work toward the lower order terms. As before, we set up two systems, one with kernel "$\widetilde{f}$" and one with kernel "$\widetilde{g}$". For the difference between the two, we use the kernel "$\widetilde{h}$". It is obvious by what has already been done, that the $\widetilde{h}_N^{i_1,\cdots,i_N}$ terms are zero. Consider an $N-1$ degree term: Fix a set of indices $\overline{i_1}, \overline{i_2}, \cdots, \overline{i_{N-1}}$.

$$\int_0^s \cdots \int_0^s \widetilde{h}_{N-1}^{\overline{i_1},\cdots,\overline{i_{N-1}}}(\sigma_1,\cdots,\sigma_{N-1}) v_{\overline{i_1}}(s-\sigma_1)\cdots v_{\overline{i_{N-1}}}(s-\sigma_{N-1})d\sigma_1\cdots d\sigma_{N-1}$$

$$+\sum_{i_N}\left\{\int_s^\infty \int_0^s \cdots \int_0^s \widetilde{f}_N^{\overline{i_1},\cdots,\overline{i_{N-1}},i_N}(\sigma_1,\cdots,\sigma_N) v_{\overline{i_1}}(s-\sigma_1)\cdots v_{\overline{i_{N-1}}}(s-\sigma_{N-1})\right.$$

$$u_{i_N}(s-\sigma_N)d\sigma_1\cdots d\sigma_N$$

$$-\int_s^\infty \int_0^s \cdots \int_0^s \widetilde{g}_N^{\overline{i_1},\cdots,\overline{i_{N-1}},i_N}(\sigma_1,\cdots,\sigma_N) v_{\overline{i_1}}(s-\sigma_1)\cdots v_{\overline{i_{N-1}}}(s-\sigma_{N-1})$$

$$\left. w_{i_N}(s-\sigma_N)d\sigma_1\cdots d\sigma_N\right\} \ .$$

The sum is over all possible $i_N$. Letting $z = u - w$ and using the polynomial argument

$$\int_0^s \cdots \int_0^s \left\{ \widetilde{h}_{N-1}^{\overline{i_1},\cdots,\overline{i_{N-1}}}(\sigma_1,\cdots,\sigma_{N-1}) + \sum_{i_N}\int_s^\infty \widetilde{f}_N^{\overline{i_1},\cdots,\overline{i_{N-1}},i_N}(\sigma_1,\cdots,\sigma_N) z(s-\sigma_N)d\sigma_N \right\}$$

$$v_{\overline{i_1}}(s-\sigma_1)\cdots v_{\overline{i_{N-1}}}(s-\sigma_{N-1})d\sigma_1\cdots d\sigma_{N-1} = 0 \ .$$

If the terms inside the parentheses are continuous,

$$\widetilde{h}_{N-1}^{\overline{i_1},\cdots,\overline{i_{N-1}}}(\sigma_1,\cdots,\sigma_{N-1}) = -\sum_{i_N}\int_s^\infty \widetilde{f}_N^{\overline{i_1},\cdots,\overline{i_{N-1}},i_N}(\sigma_1,\cdots,\sigma_N) z(s-\sigma_N)d\sigma_N \ .$$

The same limit argument used in the scalar case shows that the RHS is zero; hence,

$$\widetilde{h}_{N-1}^{\overline{i_1},\cdots,\overline{i_{N-1}}}(\sigma_1,\cdots,\sigma_{N-1}) = 0 \ .$$

Looking at what has been done so far and at the scalar case, it may be found that the kernels for the lower degree terms are zero as well. The continuity conditions are summarized below:

$$\int_s^\infty \widetilde{f}_N^{i_1,\cdots,i_N}(\sigma_1,\cdots,\sigma_N) z_{i_N}(s-\sigma_N)d\sigma_N$$

$$\int_s^\infty \int_s^\infty \widetilde{f}_N^{i_1,\cdots,i_N}(\sigma_1,\cdots,\sigma_N) z_{i_{N-1},i_N}(s-\sigma_{N-1}, s-\sigma_N)d\sigma_{N-1}d\sigma_N$$

$$\vdots$$

$$\int_s^\infty \cdots \int_s^\infty \widetilde{f}_N^{i_1,\cdots,i_N}(\sigma_1,\cdots,\sigma_N) z_{i_2,\cdots,i_N}(s-\sigma_2,\cdots, s-\sigma_N)d\sigma_2 \cdots d\sigma_N$$

are all continuous. Similarly for the lower degree terms. Sufficient conditions for continuity of these terms are similar to before; the Lebesgue Dominated Convergence Theorem is used. Also as before, the kernels need to be absolutely integrable in the polynomial case.

# APPENDIX C

# ON SHIFT DIFFERENTIABILITY OF TIME FUNCTIONS

In this appendix, we present conditions that imply the shift differentiability of a space of time functions. NSWCDD's publications format guide requires that appendices be self-contained, so we repeat here the definitions of fitted families (FFs) of time functions and other relevant definitions from the main text of this report.

**Definition C-1** *([C-1]) Let $\mathcal{L} = \mathcal{L}(\Re, E)$ be a linear space of time functions from $\Re$ into a Banach space $E$ such that any translate of a function in $\mathcal{L}$ is also a function in $\mathcal{L}$. Let $\mathcal{N} = \{\|\cdot\|_{s,t}, -\infty < s < t < \infty\}$ be a family of seminorms on $\mathcal{L}$ satisfying the following conditions:*

*(1) For $f_1, f_2 \in \mathcal{L}$, if $f_1(\tau) = f_2(\tau)$ for $s < \tau \le t$ then $\|f_1 - f_2\|_{s,t} = 0$.*
*(2) Let $\mathrm{L}_\tau$ denote shift to the left by $\tau$. For all $f \in \mathcal{L}$, $\|\mathrm{L}_\tau f\|_{s-\tau,t-\tau} = \|f\|_{s,t}$.*
*(3) Let $r < s < t$. Then for all $f \in \mathcal{L}$, $\|f\|_{s,t} \le \|f\|_{r,t}$.*
*(4) Let $r < s < t$. Then for all $f \in \mathcal{L}$, $\|f\|_{r,t} \le \|f\|_{r,s} + \|f\|_{s,t}$.*
*(5) There exists $0 < \alpha \le \infty$ and $K \ge 1$ such that if $0 < t - r \le \alpha$ and $r < s < t$, then for all $f \in \mathcal{L}$, $\|f\|_{r,s} \le K \|f\|_{r,t}$.*

*The pair $(\mathcal{L}, \mathcal{N})$ is called a fitted family of seminorms on $\mathcal{L}$. The normed linear space formed from equivalence classes of functions in $\mathcal{L}$ with norm $\|\cdot\|_{s,t}$ is denoted $A_{s,t}$. The elements of $A_{s,t}$ are the equivalence classes determined by: $f \sim g$, $f, g \in \mathcal{L}$ if and only if $\|f - g\|_{s,t} = 0$. They are denoted $u_{s,t}$, $y_{s,t}$, etc. The set $\{A_{s,t}\}$, $-\infty < s < t < \infty$, is the FF of normed linear spaces given by $(\mathcal{L}, \mathcal{N})$.*

A FF $(\mathcal{L}, \mathcal{N})$ and $\{A_{s,t}\}$, $-\infty < s < t < \infty$, can be augmented to include $\|\cdot\|_{-\infty,t}$ by taking the limit $s \to -\infty$, since by (3) of Definition C-1 $\|f\|_{s,t}$ is monotone nondecreasing as $s \to -\infty$ with $t$ fixed. Let $\mathcal{L}_0 = \{f \in \mathcal{L} |\ \lim_{s\to-\infty} \|f\|_{s,t} < \infty,\ t \in \Re\}$. For $f \in \mathcal{L}_0$, define

$$\|f\|_t \triangleq \lim_{s\to-\infty} \|f\|_{s,t} = \|f\|_{-\infty,t}\,. \tag{C-1}$$

Let $\mathcal{L}_{00} = \{f \in \mathcal{L}_0 |\ \sup_t \|f\|_t < \infty\}$. For $f \in \mathcal{L}_{00}$ define

$$\|f\|_{s,\infty} \triangleq \sup_{t>s} \|f\|_{s,t}\ ;\ -\infty \le s\,. \tag{C-2}$$

For $f \in \mathcal{L}_{00}$ we put

$$\|f\| \triangleq \sup_{t\in\Re} \|f\|_t = \|f\|_{-\infty,\infty}\,. \tag{C-3}$$

The normed linear space consisting of equivalence classes of functions in $\mathcal{L}_{00}$ with the norm (C-3) is called the *bounding space* $A$ for the family $\{A_{s,t}\}$.

**Definition C-2** *The space of time functions $U$ is shift differentiable if for all $u \in U$*

$$\mathrm{L}_h u = u + h \cdot d + e(h) \tag{C-4}$$

*where $d, e(h) \in U$ and $\lim_{h\to 0} \|e(h)\|_t / h = 0$ for all $t \in \Re$. $d$ is the shift derivative of $u$.*

Shift continuity, as it is stated in [C-2], is the space $U$ is *shift continuous* if $\lim_{h\to0} \|u - \mathrm{L}_h u\|_{s,t} = 0$ for all $u$ for which the norm is defined and for all $-\infty \le s < t < \infty$. From Definition C-2, we see shift differentiability implies shift continuity.

Before going to our result, we start out by showing that differentiability or continuity is not a requirement for shift differentiability by these two examples:

**Example C-3** *Let the basic function space $\mathcal{L}(\Re, \Re)$ be the set of real-valued functions on $\Re$ that are bounded and piecewise continuous (continuous except possibly at a finite number of points on any finite interval and with finite left- and right-hand limits at these points.) The family of norms is given by the $L_2$ norm $\|f\|_{s,t} = \int_s^t |f(\tau)|^2 d\tau$, $f \in \mathcal{L}(\Re, \Re)$. The (usually input) space $U$ is either the bounding space $A$ for this family or a shift-invariant subset thereof. Define $u$ and $d$ by*

$$u(t) = \begin{cases} 0, & t \leq -1; t > 11 \\ t+1, & -1 < t \leq 0 \\ 1, & 0 < t \leq 10 \\ -t+11, & 10 < t \leq 11 \end{cases}$$

$$d(t) = \begin{cases} 0, & t \leq -1 \\ 1, & -1 < t \leq 0 \\ 0, & 0 < t \leq 10 \\ -1, & 10 < t \leq 11 \\ 0, & 11 < t \end{cases} .$$

*We observe that $u$ and $d$ are in $L_2$. We note that $u$ is not differentiable; however, we show below that $u$ is shift differentiable. We have that (with $0 < h < 10$)*

$$\mathrm{L}_h u(t) = \begin{cases} 0, & t \leq -1-h \\ t+(1+h), & -1-h < t \leq -h \\ 1, & -h < t \leq 10-h \\ -t+(11-h), & 10-h < t \leq 11-h \\ 0, & 11-h < t \end{cases} .$$

*Then we calculate*

$$(\mathrm{L}_h u - u - h \cdot d)(t) = \begin{cases} 0, & t \leq -1-h \\ t+(1+h), & -1-h < t \leq -1 \\ 0, & -1 < t \leq -h \\ -(t+h), & -h < t \leq 0 \\ 0, & 0 < t \leq 10-h \\ -t+(10-h), & 10-h < t \leq 10 \\ 0, & 10 < t \leq 11-h \\ t-(11-h), & 11-h < t \leq 11 \\ 0, & 11 < t \end{cases} .$$

*Noting that $e = (\mathrm{L}_h u - u - h \cdot d)$. We first calculate*

$$\|\mathrm{L}_h u - u - h \cdot d\|^2 = 2 \cdot 2 \cdot \int_{-h}^{0} t^2 dt = \frac{4 \cdot h^3}{3}$$

*then*

$$\lim_{h\to 0}\frac{\|e(h)\|}{h}=\lim_{h\to 0}\frac{\sqrt{\dfrac{4\cdot h^3}{3}}}{h}=\frac{2}{\sqrt{3}}\lim_{h\to 0}\frac{h^{3/2}}{h}=0 \quad .$$

**Example C-4** *Let the basic function space $\mathcal{L}(\Re,\Re)$ for inputs be the set of real-valued functions on $\Re$ that are bounded and piecewise continuous (continuous except possibly at a finite number of points on any finite interval and with finite left- and right-hand limits at these points). The family of norms is given by $L_1$ norm $\|f\|_{s,t}=\int_s^t |f(\tau)|d\tau$, $f\in\mathcal{L}(\Re,\Re)$. The (usual input) space $U$ is either the bounding space $A$ for this family or a shift-invariant subset thereof. This time let $u$ be the boxcar step and $d$ be the addition of delta functions, respectively,*

$$u(t)=\begin{cases}1, & 0<t\le 10\\ 0, & \textit{elsewhere}\end{cases}$$

$$\delta(t)=\begin{cases}0, & t\neq 0\\ \textit{undefined}, & t=0\end{cases}$$

*with*

$$\int_{-\infty}^{\infty}\delta(\tau)\,d\tau=1 \quad .$$

*(The usual symbol for the delta function is $\delta$.)*

$$d(t)=\delta(t)-\delta(t-10)$$

*We note that $u$ is not continuous; however, we show below that $u$ is shift differentiable. We have that (with $h>0$)*

$$\mathrm{L}_h u(t)=\begin{cases}1, & -h<t\le 10-h\\ 0, & \textit{elsewhere}\end{cases} \quad .$$

*Then we calculate $e=(\mathrm{L}_h u-u-h\cdot d)$. First define $g$ by*

$$[g(h)](t)=(L_h u-u)(t)=\begin{cases}0, & t\le -h\\ 1, & -h<t\le 0\\ 0, & 0<t\le 10-h\\ -1, & 10-h<t\le 10\\ 0, & 10<t\end{cases} \quad .$$

*Then*

$$[e(h)](t)=[g(h)](t)-hd(t)=[g(h)](t)-h\delta(t)+h\delta(t-10) \quad .$$

*Looking at the essential property $e$ must have,*

$$\lim_{h\to 0}\frac{\|[e(h)](t)\|}{h}=\lim_{h\to 0}\int_{-\infty}^{\infty}\frac{|[e(h)](\tau)|}{h}d\tau$$

$$=\lim_{h\to 0}\frac{1}{h}\left\{\int_{-h}^{0}1\cdot d\tau-\int_{10-h}^{10}1\cdot d\tau-\int_{-\infty}^{\infty}h\cdot d(\tau)\,d\tau\right\}$$

$$= \lim_{h\to 0} \frac{1}{h} \left\{ \int_{-h}^{0} 1 \cdot d\tau - h \cdot \int_{-\infty}^{\infty} \delta(\tau)\, d\tau \right\} - \lim_{h\to 0} \frac{1}{h} \left\{ \int_{10-h}^{10} 1 \cdot d\tau - h \cdot \int_{-\infty}^{\infty} \delta(\tau - 10)\, d\tau \right\}$$

$$= 2 \cdot \lim_{h\to 0} \left[ \frac{1}{h} \{h - h\} \right] = 2 \cdot \lim_{h\to 0} \left[ \frac{1}{h} \cdot 0 \right] = 2 \cdot \lim_{h\to 0} [0] = 0 \quad .$$

The next proposition gives sufficient conditions for shift differentiability.

**Proposition C-5** *Let $\mathcal{L} = \mathcal{L}(\Re, E)$ be a linear space of time functions. Use the $L_p$ norm. For $f \in \mathcal{L}$, let $f'(t)$ be the derivative, that is,*

$$f'(t) = \lim_{h\to 0} \frac{f(t+h) - f(t)}{h} \tag{C-5}$$

*for all $t \in \Re$. Define*

$$e(t,h) \equiv f(t+h) - f(t) - hf'(t)$$

*for all $t$, $h \in \Re$. Assume $e(t,\cdot)$ is measurable. Assume there exists $\epsilon(t) \in L_1(-\infty, T]$ such that*

$$\left\| \frac{e(t,h)}{h} \right\| \le \epsilon(t) \tag{C-6}$$

*for all $t$, $h \in \Re$. Then*

$$\mathrm{L}_h f(t) = f(t+h) = f(t) + hf'(t) + e(t,h) \tag{C-7}$$

*with*

$$\lim_{h\to 0} \frac{\|e(h)\|_T}{h} = \lim_{h\to 0} \frac{1}{h} \sqrt[p]{\int_{-\infty}^{T} \|e(t,h)\|^p\, dt} = \lim_{h\to 0} \sqrt[p]{\int_{-\infty}^{T} \left\| \frac{e(t,h)}{h} \right\|^p dt}$$

*using (C-6) and the Lebesgue Dominated Convergence Theorem,*

$$= \sqrt[p]{\int_{-\infty}^{T} \lim_{h\to 0} \left\| \frac{e(t,h)}{h} \right\|^p dt} = \sqrt[p]{\int_{-\infty}^{T} \left\| \lim_{h\to 0} \frac{e(t,h)}{h} \right\|^p dt} = 0 \ .$$

*So, $f'$ is the shift derivative of $f$. If the limit in (C-5) and the bound in (C-6) are almost everywhere, then the result is in terms of the eqivalence classes $u \in U$,*

$$\mathrm{L}_h u = u + hd + e(h) \ . \tag{C-8}$$

# APPENDIX D

# ON SHIFT DIFFERENTIABILITY OF SYSTEM TRAJECTORIES

In this appendix, sufficient conditions for a system trajectory to be shift-differentiable are considered. That the system be time-invariant is a straightforward sufficient condition. For time-varying systems, our conditions are based on a hypothesis SG, defined below, and the semigroup it defines being strongly continuous. NSWCDD's publications format guide requires that appendices be self-contained, so we recover here several definitions from the main text of this report.

An input-output system is denoted $(Y, F, U)$, where $F$ is a mapping from an input space $U$ to an output space $Y$, and where $U$ and $Y$ are translation-invariant spaces of vector-valued time functions. The input and output space metrics are set up by seminorms referred to as fitted families (FFs) of seminorms. The notation $\|u_{s,t}\|_{s,t}$ indicates the norm of the input $u$ over the interval of time $(s, t]$. $U_{s,t}$ is the space of inputs over the same interval. FFs were initially described in [D-1].

**Definition D-1** *[D-1] Let $\mathcal{L} = \mathcal{L}(\Re, E)$ be a linear space of time functions from $\Re$ into a Banach space $E$ such that any translate of a function in $\mathcal{L}$ is also a function in $\mathcal{L}$. Let $\mathcal{N} = \{\|\cdot\|_{s,t}, -\infty < s < t < \infty\}$ be a family of seminorms on $\mathcal{L}$ satisfying the following conditions:*

*(1) For $f_1, f_2 \in \mathcal{L}$, if $f_1(\tau) = f_2(\tau)$ for $s < \tau \leq t$ then $\|f_1 - f_2\|_{s,t} = 0$.*

*(2) Let $\mathrm{L}_\tau$ denote shift to the left by $\tau$. For all $f \in \mathcal{L}$, $\|\mathrm{L}_\tau f\|_{s-\tau,t-\tau} = \|f\|_{s,t}$.*

*(3) Let $r < s < t$. Then for all $f \in \mathcal{L}$, $\|f\|_{s,t} \leq \|f\|_{r,t}$.*

*(4) Let $r < s < t$. Then for all $f \in \mathcal{L}$, $\|f\|_{r,t} \leq \|f\|_{r,s} + \|f\|_{s,t}$.*

*(5) There exists $0 < \alpha \leq \infty$ and $K \geq 1$ such that if $0 < t - r \leq \alpha$ and $r < s < t$, then for all $f \in \mathcal{L}$, $\|f\|_{r,s} \leq K \|f\|_{r,t}$.*

*The pair $(\mathcal{L}, \mathcal{N})$ is called a FF of seminorms on $\mathcal{L}$. The normed linear space formed from equivalence classes of functions in $\mathcal{L}$ with norm $\|\cdot\|_{s,t}$ is denoted $A_{s,t}$. The elements of $A_{s,t}$ are the equivalence classes determined by: $f \sim g$, $f, g \in \mathcal{L}$ if and only if $\|f - g\|_{s,t} = 0$. They are denoted $u_{s,t}$, $y_{s,t}$, etc. The set $\{A_{s,t}\}$, $-\infty < s < t < \infty$, is the FF of normed linear spaces given by $(\mathcal{L}, \mathcal{N})$.*

We also define $\|u\|_{-\infty,t} = \lim_{s\to-\infty} \|u\|_{s,t} = \|u\|_t$. Also, $\|f\|_{s,\infty} = \sup_{t>s} \|f\|_{s,t}$, and $\|f\| = \sup_t \|f\|_t$ . A mapping $F : U \to Y$ is called a *global input-output mapping.*

**Definition D-2** *Let $(Y, F, U)$ be an input-output system. $F$ is a causal mapping and $(Y, F, U)$ is a causal system if and only if for all $t$ and for all $u, v \in U$ such that $\|u - v\|_t = 0$; it follows that $\|F(u) - F(v)\|_t = 0$.*

If $F$ satisfies this definition, it determines a mapping from $U_t$ into $B_t$, denoted $\widetilde{F}_t$, that satisfies $\left\|\widetilde{F}_t u_t - (Fu)_t\right\|_t = 0$. We call $\widetilde{F}_t$ a *truncated input-output mapping*, and define the *centered truncated input-output mapping* $\underset{\sim}{F}_t : U_0 \to Y_0$ by $\underset{\sim}{F}_t(u_0) \overset{\Delta}{=} \mathrm{L}_t \widetilde{F}_t \mathrm{R}_t(u_0)$, where $\mathrm{R}_t \overset{\Delta}{=} \mathrm{L}_{-t}$ is the right-shift by $t$. The norms we use here for input-output mappings $F$, $\widetilde{F}_t$, and $\underset{\sim}{F}_t$ are the $N$-power norms, denoted $[\![ \cdot ]\!]$.

Let $\Phi$ be a mapping from a normed linear space $X$ into a normed linear space $Y$. For any nonnegative integer $N$, the $N$-power norm for $\Phi$ is given by

$$[\![\Phi]\!] \overset{\Delta}{=} \sup_{x\in X} \frac{\|\Phi(x)\|}{1 + \|x\|^N} \qquad \text{(D-1)}$$

when the right side exists. We say $\Phi$ is *bounded (in $N-$power norm)* if $\|\Phi\| < \infty$. Let $\mathcal{F}_N(X,Y)$ be the normed linear space of all mappings $\Phi : X \to Y$ with $\|\Phi\| < \infty$, and $\mathcal{C}_N(X,Y)$ be the normed linear subspace of $\mathcal{F}_N(X,Y)$ of all continuous $\Phi$.

**Definition D-3** *A system trajectory is a shift differentiable system trajectory $t \to F_t$ if*

$$\mathrm{L}_h F \mathrm{R}_h = F + h \cdot G + H(h) \tag{D-2}$$

*where $G \in \mathcal{C}_N(U,Y)$ and $H(h) \in \mathcal{C}_N(U,Y)$ and such that*

$$\lim_{h \to 0} \frac{[]H(h)_t[]}{h} = 0 \text{ for all } t \in \Re \cdot \tag{D-3}$$

We define the truncation mapping $\pi_t$ on the space of systems $\mathcal{C}_N(X,Y)$ by $\pi_t F = F_t$.

**Definition D-4** *We say $\mathbf{F} = \{(Y,F,U)\}$ satisfies condition (SG) if for all $F^1, F^2 \in \mathbf{F}$, $\pi_t F^1 = \pi_s F^2$ for some $t, s \in \Re$ implies $\pi_{t+a} F^1 = \pi_{s+a} F^2$ for all $a \geq 0$.*

Now, define an operator $\theta(s,t)$ for $s,t \in \Re$, $s \leq t$, with domain $\mathbf{F}_s$ and counter-domain $\mathbf{F}_t$ by

$$\theta(s,t) F_s = \pi_t \circ \pi_s^{-1} F_s \; .$$

If condition (SG) is satisfied, $\theta(s,t)$ is a well-defined operator that is linear on its domain. In fact, $\theta$ only depends on $(t-s)$. The family of operators $\{\theta(\tau); \tau \geq 0\}$ is a one-parameter semigroup of operators, linear on the domain specified.

**Lemma D-5** *[D-2, Lemma 5.5.1] Let $\mathbf{F} = \{(Y,F,U)\}$ be as described in the preceding paragraph. Let $\{\theta(\tau), \tau > 0\}$ be the semigroup defined by condition (SG). Assume that $\{\theta(\tau), \tau > 0\}$ is a strongly continuous semigroup of bounded linear transformations with infinitesimal generator $A$. For $F \in \mathbf{F}$, if $F_t \in Domain(A)$ for all $t$ and*

$$\lim_{h \to 0} \frac{\theta(h) F_t - F_t}{h} - A F_t = 0 \tag{D-4}$$

*uniformly in $t$, then the system trajectory $t \to F_t$ is shift differentiable.*

*Proof*: Consider the system trajectory $t \to F_t$ of the input-output system $(Y,F,U)$. Fix $h > 0$ and consider the time function $t \to \theta(h) F_t$ or $t \to F_{t+h}$, which is a time function in $\mathcal{L}(\Re, \mathcal{C}_N(U_0, Y_0))$. The time function $t \to F_{t+h}$ is actually defined by the causal input-output system $(Y, \mathrm{L}_h F \mathrm{R}_h, U)$. By the uniform convergence in (D-4), the sequence

$$\left\{ \frac{L_{(1/n)} F R_{(1/n)} - F}{(1/n)} \right\}_{n=1}^{+\infty}$$

is Cauchy in the metric of $\mathcal{C}_N(X,Y)$. Let this sequence converge to $G \in \mathcal{C}_N(X,Y)$. (The completeness of $\mathcal{C}_N(X,Y)$ is considered in Lemma A.1 of [D-2].) Note that for all $t \in \Re$, $A \cdot F_t = G_t$. For each $h \in \Re$, define $H(h) \in \mathcal{C}_N(U_0, Y_0)$ by

$$H(h) = L_h F R_h - F - hG$$

It is seen that $\lim_{h\to 0} H(h)/h = 0$; therefore, the system trajectory $t \to F_t$ is shift differentiable.
□

Lemma D-5 gives (D-3) is uniform in $t$.

## DISTRIBUTION

Copies

**DOD ACTIVITIES (CONUS)**

DEFENSE TECHNICAL 2
INFORMATION CENTER
8725 JOHN J KINGMAN ROAD
SUITE 0944
FORT BELVOIR VA 22060-6218

DEBORAH D HUBBARD 1
OUSD(AT&L) / IC / CWP
3070 DEFENSE PENTAGON
PENTAGON SUITE 2E173A
WASHINGTON DC 20301-3070

**NON-DOD ACTIVITIES (CONUS)**

DENNIS S BERNSTEIN 1
PROFESSOR
AEROSPACE ENGINEERING DEPARTMENT
1320 BEAL ST
THE UNIVERSITY OF MICHIGAN
ANN ARBOR MI 48109

DOCUMENT CENTER 1
THE CNA CORPORATION
4825 MARK CENTER DRIVE
ALEXANDRIA VA 22311-1850

JOHN CHIN 4
GOVERNMENT DOCUMENTS SECTION
101 INDEPENDENCE AVENUE SE
LIBRARY OF CONGRESS
WASHINGTON DC 20540-4172

GLENN K HEITMAN 1
4830 NW 43RD STREET
GAINESVILLE FL 32606

PIERRE T KABAMBA 1
PROFESSOR
AEROSPACE ENGINEERING DEPARTMENT
1320 BEAL ST
THE UNIVERSITY OF MICHIGAN
ANN ARBOR MI 48109

Copies

LAWRENCE C NG 1
LOCKHEED MARTIN SPACE SYSTEMS
DEPARTMENT DACS BUILDING 150
1111 LOCKHEED MARTIN WAY
SUNNYVALE CA 94089

FRANK RIEFLER 1
LOCKHEED MARTIN CORPORATION
MS2 137-205
199 BORTON LANDING ROAD
MOORESTOWN NJ 08057-3075

KHANH D PHAM 1
DYNAMICS AND CONTROLS
AIR FORCE RESEARCH LABORATORY
SPACE VEHICLES DIRECTORATE
AFRL/VSSV
3550 ABERDEEN AVE SE
KIRTLAND AFB NM 87117

JASON SPEYER 1
PROFESSOR
MECHANICAL & AEROSPACE
ENGINEERING DEPARTMENT
UNIVERSITY OF CALIFORNIA LOS ANGELES
420 WESTWOOD PLAZA 38-137O
LOS ANGELES CA 90095-1597

LESZEK J SCZANIECKI 1
LOCKHEED MARTIN CORPORATION
MS2 137-205
199 BORTON LANDING ROAD
MOORESTOWN NJ 08057-3075

JOHN STETSON 1
LOCKHEED MARTIN CORPORATION
MS2 137-205
199 BORTON LANDING ROAD
MOORESTOWN NJ 08057-3075

ANDREW VOGT 1
ASSOCIATE PROFESSOR
DEPARTMENT OF MATHEMATICS
AND STATISTICS
BOX 571233
ST. MARY'S HALL 338A
GEORGETOWN UNIVERSITY
WASHINGTON DC 20057-1223

**DISTRIBUTION (Continued)**

| | Copies |
|---|---|
| CHRISTOFER ZAFFRAM<br>MCCDC STUDIES AND ANALYSIS<br>BRANCH<br>3300 RUSSELL ROAD<br>QUANTICO VA 22134 | 1 |

**NON-DOD ACTIVITIES (EX-CONUS)**

| | |
|---|---|
| DRUMI BAINOV<br>PROFESSOR<br>MEDICAL UNIVERSITY OF SOFIA<br>PO BOX 45<br>SOFIA 1504 BULGARIA | 1 |
| JOACHIM ROSENTHAL<br>PROFESSOR<br>DEPARTMENT OF MATHEMATICS<br>UNIVERSITY OF ZURICH<br>WINTERTHURERSTRASSE 190<br>CH-8057 ZURICH | 1 |

**INTERNAL**

| | |
|---|---|
| K94 (LAMBERTSON) | 1 |
| Q23 (BILLARD) | 1 |
| Q23 (GRAY) | 1 |
| Q31 (BREAUX) | 1 |
| Q40 (YOST) | 1 |
| Q45 (ROGINSKY) | 1 |
| W11 (ORMSBY) | 1 |
| W30 (DIDONATO) | 1 |
| W30 (SKAMANGAS) | 1 |
| W31 (BECK) | 1 |
| W31 (CARR) | 1 |
| W31 (LAWTON) | 1 |
| W31 (HABGOOD) | 1 |
| W31 (MCDEVITT) | 1 |
| W31 (SERAKOS) | 5 |
| W31 (YOUSSEF) | 1 |
| W41 (MCCABE) | 1 |
| W42 (POLEY) | 1 |
| W43 (AFRICA) | 1 |
| Z31 (TECHNICAL LIBRARY) | 3 |

NAVSEA

WARFARE CENTERS

DAHLGREN